\documentclass[10pt]{article}
\usepackage{amssymb,amsmath,amsfonts,amsthm}
\usepackage[dvips]{graphicx}
    \newtheorem{theorem}{Theorem}[section]
    
    \newtheorem{lemma}[theorem]{Lemma}
    \newtheorem{proposition}[theorem]{Proposition}
    \newtheorem{remark}[theorem]{Remark}
    \newcommand{\U}{\mathbf{u}}
    \newcommand{\V}{\mathbf{v}}
    \newcommand{\W}{\mathbf{w}}
    \newcommand{\Z}{\mathbf{z}}
\begin{document}
\title{Perturbation from symmetry and multiplicity of solutions for
 strongly indefinite elliptic systems}
\author{Cristina Tarsi \thanks{ e-mail:
Cristina.Tarsi@mat.unimi.it.}\\
Dipartimento di Matematica, Universit\`{a} degli Studi, \\
I-20133 Milano, Italy}
\date{}
\maketitle {\footnotesize {\bf Abstract} \noindent We consider the
following elliptic system:
\[
\left\{
\begin{array}{ll}
 -\Delta u= |v|^{p-1} v + h(x) %
& \;\;\;\;\;x\in \Omega \\
 -\Delta v= |u|^{q-1} u + k(x) %
& \;\;\;\;\;x\in \Omega \\
u=v=0 & \;\;\;\;\; x\in \partial \Omega
\end{array}
\right.
\]
where $\Omega \subset \mathbb{R} ^N$, $N\geq 3$ is a smooth
bounded domain. If $h(x)\equiv k(x)\equiv 0$ the system presents a
natural $\mathbb{Z}_2$ symmetry, which guarantees the existence of
infinitely many solutions. In this paper we show that the
multiplicity structure can be maintained if $(p,q)$ lies below a
suitable curve in $\mathbb{R}^2$. }
\section{Introduction}
There has been recently an active research in the study of
semilinear elliptic systems: see for example \cite{dF} for a
survey on the argument. Such systems are called \emph{variational}
if solutions can be viewed as critical points of an associated
functional defined on a suitable function space. Restricting our
attention to second order elliptic systems with two unknowns,
whose principal part is given by the differential operator
$-\Delta$, we consider systems of the form
\begin{equation}\label{genericsyst}
\left\{
\begin{array}{ll}
 -\Delta u= f(x,u,v) %
& \;\;\;\;\;x\in \Omega \\
 -\Delta v= g(x,u,v) %
& \;\;\;\;\;x\in \Omega \\
u=v=0 & \;\;\;\;\; x\in \partial \Omega
\end{array}
\right.
\end{equation}
where $\Omega$ is a smooth bounded domain in $\mathbb{R} ^N$,
$N\geq 3$.

We say that \eqref{genericsyst} is a \emph{potential system} if
there exists a function $F:\overline{\Omega}\times \mathbb{R}
\times \mathbb{R} \rightarrow \mathbb{R}$ of class $C^1$ such that
\[
\frac{\partial F}{\partial u}=f, \;\;\;\;\; %
-\frac{\partial F}{\partial v}=g,
\]
that is,
\[
\left\{
\begin{array}{ll}
 -\Delta u= \partial_u F %
& \;\;\;\;\;x\in \Omega \\
 +\Delta v= \partial_v F %
& \;\;\;\;\;x\in \Omega \\
u=v=0 & \;\;\;\;\; x\in \partial \Omega.
\end{array}
\right.
\]
These are the Euler-Lagrange equations of the functional
\[
\Phi(u,v)=\frac12\int_{\Omega}|\nabla u|^2-
\frac12\int_{\Omega}|\nabla v|^2-\int_{\Omega} F(x,u,v)
\]
whose critical points are the weak solutions of equations
\eqref{genericsyst}. If $F$ satisfies suitable growth conditions,
this functional is well defined in the cartesian product
$E=H_0^1(\Omega)\times H_0^1(\Omega)$, by virtue of the Sobolev
embedding theorem; note that $\Phi$ has a strongly indefinite
quadratic part. Systems of this type have been studied, for
example, in \cite{BR}, \cite{CM}, \cite{dFF}; recently existence
and multiplicity results have been obtained also for indefinite
systems with critical growth (see e.g. \cite{dFD}, \cite{CDH-L}).

We say that \eqref{genericsyst} is a \emph{Hamiltonian system} if
there exists a function $H:\overline{\Omega}\times \mathbb{R}
\times \mathbb{R} \rightarrow \mathbb{R}$ of class $C^1$ such that
\[
\frac{\partial H}{\partial v}=f, \;\;\;\;\; %
\frac{\partial H}{\partial u}=g,
\]
that is,
\[
\left\{
\begin{array}{ll}
 -\Delta u= \partial_v H %
& \;\;\;\;\;x\in \Omega \\
 -\Delta v= \partial_u H %
& \;\;\;\;\;x\in \Omega \\
u=v=0 & \;\;\;\;\; x\in \partial \Omega.
\end{array}
\right.
\]
By analogy with the scalar case one would guess that the
subcritical case occurs if the growths of $H$ with respect to $u$
and $v$ are both less than $2^{*}=(N+2)/(N-2)$: in this case one
could search the weak solutions of the Hamiltonian system as
critical points of the functional
\[
\Phi(u,v) = \int_{\Omega}\nabla u \nabla v -\int_{\Omega} H(x,u,v)
\]
which is well defined on $E=H_0^1(\Omega)\times H_0^1(\Omega)$.
Nevertheless, the coupling now also occurs in the quadratic part
of $\Phi$, and therefore is much stronger than in the potential
case. An immediate consequence is that this approach is too
restrictive: there is no longer one appropriate choice of function
spaces, and the notion of criticality have to take into
consideration the fact that the system is coupled. In \cite{CdFM},
\cite{dFF} and \cite{HvdV} appeared the notion of \emph{Critical
Hyperbola}, which replaces the notion of critical exponent of the
scalar case when $N\geq3$,
\begin{equation}\label{crithyp}
\frac 1{p+1}+\frac 1{q+1}=1-\frac2{N}
\end{equation}
that is associated to Hamiltonian system when $\partial_v H$ grows
like $v^p$ as $v\rightarrow +\infty$ and $\partial_u H$ grows like
$u^q$ as $u\rightarrow +\infty$, and the dependance on the other
variables is of some lower orders. It is known that for any point
$(p,q)\in \mathbb{R}^2$ below the critical hyperbola the
Hamiltonian system has a nontrivial solution (see \cite{CdFM},
\cite{dFF}, \cite{HvdV}, \cite{FM} and \cite{dFR}), whereas for
points $(p,q)$ on the critical hyperbola one finds the typical
problems of non-compactness and non-existence of solutions (see
\cite{vdV} and \cite{M}).

If $F(x,u,v)$ or $H(x,u,v)$ is even in $(u,v)$, the potential,
respectively Hamiltonian, system possesses a natural $\mathbb{Z}
_2$-symmetry: by analogy with the scalar case, one would expect
the existence of infinitely many solutions. In the scalar case,
the standard variational method for dealing with even equations is
based on the symmetric version of the Mountain Pass Theorem of
Ambrosetti-Rabinowitz (see \cite{St}); this theorem is no longer
applicable in the case of an elliptic system, since the functional
associated is strongly indefinite. Nevertheless, by means of a
Galerkin type approximation, one can reduce the strongly
indefinite functional to a semidefinite situation (see \cite{BW},
\cite{BC}, \cite{D}, \cite{BdF}, \cite{dFD} and others). A
different approach to the problem of symmetric indefinite
functional was given by Angenent and van der Vorst in \cite{AvdV},
who applied Floer's version of Morse theory to Hamiltonian
elliptic systems, in the spirit of \cite{BL}; see also
\cite{AvdV2}.

As in the scalar case, one could ask if the multiplicity structure
can be maintained by adding a perturbation term of lower order.
This problem has been extensively investigated in the case of a
single equation: a partial answer was independently obtained by
Struwe \cite{St1}, Bahri-Berestycki \cite{BB}, Rabinowitz
\cite{Ra}, Bahri-Lions \cite{BL}, who showed in important works
that the multiplicity structure can be maintained restricting the
growth range of the nonlinearity with suitable bounds depending on
$N$.
\\The problem of perturbation from
symmetry of elliptic systems have been treated, to our knowledge,
only by Clapp, Ding and Hernandez-Linares in \cite{CDH-L}. Here
the authors obtain a multiplicity result only for perturbed
symmetric potential systems, where the perturbation terms can
depend also on the unknowns $(u,v)$ (with suitable limitations on
the growth in $u$, $v$). The proof, as mentioned before, is based
on a Galerkin type approximation, which reduces the study of the
strongly indefinite functional associated to the potential system
to a semidefinite situation, thus allowing the use of the Morse
theory methods as in \cite{BL}.

The aim of this paper is to obtain a multiplicity result for
Hamiltonian systems with perturbed symmetries of the type:
\begin{equation}
\label{pertsystem} \left\{
\begin{array}{ll}
 -\Delta u= |v|^{p-1} v + h(x) %
& \;\;\;\;\;x\in \Omega \\
 -\Delta v= |u|^{q-1} u + k(x) %
& \;\;\;\;\;x\in \Omega \\
u=v=0 & \;\;\;\;\; x\in \partial \Omega
\end{array}
\right.
\end{equation}
where $\Omega \subset \mathbb{R} ^N$, $N\geq 3$, is a smooth
bounded domain. In this case, as observed, there is no longer a
one appropriate choice of the function spaces: in \cite{HvdV} and
\cite{dFF} the authors propose the use of Sobolev spaces of
fractional order, obtaining the critical hyperbola, whereas in
\cite{dFdOR} the authors choose a different approach, based on an
Orlicz space setting, which yields to the same result when the
hypotheses overlap. In this paper we follow the idea of de
Figueiredo-Felmer \cite{dFF} and van der Vorst \cite{HvdV},
defining the variational setting of \eqref{pertsystem} on a
cartesian product of suitable fractional Sobolev spaces: roughly
speaking, these spaces, denoted by $H^s(\Omega)$, $s>0$, consist
of functions whose derivative of order $s$ is in $L^2(\Omega)$
(they can be defined by means of interpolation or Fourier
expansion). Therefore, even if we reduce to a semidefinite
situation by means of Galerkin type approximation, the classical
Morse theory methods, as in \cite{BL} and in \cite{CDH-L}, are not
applicable. For this reason our approach to the problem of
perturbation from symmetry follows the first one proposed by
Struwe \cite{St1}, Bahri-Berestycki \cite{BB} and Rabinowitz
\cite{Ra}. This yields the following theorem, which is our main
result:
\begin{theorem}\label{mainresult}
Let $\Omega$ be a smooth bounded domain in $\mathbb{R} ^N$, $N\geq
3$, and let $p,q>1$ satisfying the following conditions
\begin{eqnarray}\label{condTh}
\nonumber \frac1{p+1}+\frac1{q+1}+\frac{p+1}{p(q+1)} &>&\frac{2N-2}{N} \hspace{10 pt}%
\textrm{if } q\geq p \\
\\
\nonumber \frac1{p+1}+\frac1{q+1}+\frac{q+1}{q(p+1)} &>&\frac{2N-2}{N} \hspace{10 pt}%
\textrm{if } q\leq p.
\end{eqnarray}
Then, for any $h,k \in L^2(\Omega)$ problem \eqref{pertsystem} has
infinitely many solutions.
\end{theorem}
Conditions \eqref{condTh} define a region in the $(p,q)$ plane
which is strictly contained in the subcritical one, as shown in
Fig. 1.
\begin{figure}[h]
\begin{center}
\includegraphics[width=0.6\textwidth]{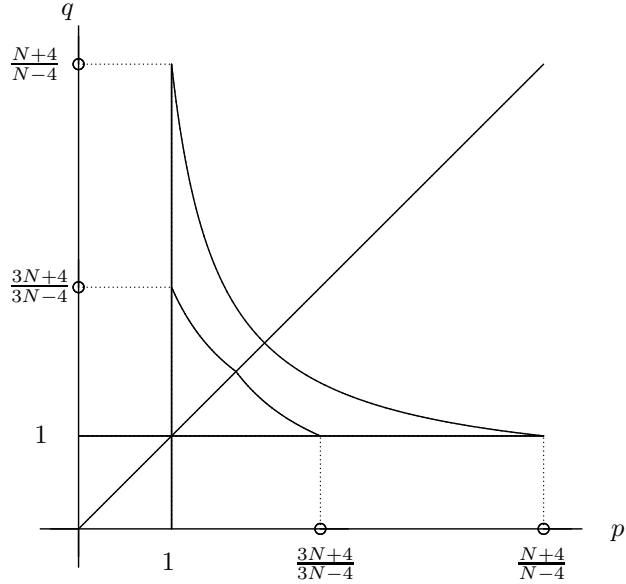}
\put(10,13){$p$} \put(-26,0){$\frac{N+4}{N-4}$}
\put(-110,0){$\frac{3N+4}{3N-4}$} \put(-160,0){1} \put(-208,48){1}
\put(-218,105){$\frac{3N+4}{3N-4}$}
\put(-218,189){$\frac{N+4}{N-4}$} \put(-198,210){$q$}
\caption{Critical Hyperbola and Theorem \eqref{mainresult},
$N>4$.}
\end{center}
\end{figure}

The paper is organized as follows. In Section 2 we describe the
variational formulation of the perturbed system
\eqref{pertsystem}. In Section 3 we consider the symmetric system
which arises from \eqref{pertsystem} when $h(x)\equiv k(x)\equiv
0$, exhibiting the unbounded sequence of critical values of the
functional associated to the symmetric problem. In Section 4 we
define a suitable modified functional $J$ associated to the
perturbed problem, in the spirit of \cite{Ra2}, whose critical
points are solutions of \eqref{pertsystem}. In Section 5 we
construct minimax sequences strictly related to the existence of
critical points of $J$, by means of Galerkin type approximation.
In Section 6  and 7, finally, we prove Theorem \ref{mainresult}
comparing upper and lower bounds of the minimax sequences
constructed before.
\begin{remark}
Another possible approach to the perturbed Hamiltonian system
\eqref{pertsystem} could be trying to apply Floer's version of
Morse theory as done by Angenent and van der Vorst in \cite{AvdV},
in the spirit of \cite{BL}. However, in the scalar case the proof
deeply depends on the relation between the Morse index of a
critical point and the number of non positive eigenvalues of the
operator $-\Delta +V(x)$ (see \cite{BL} or \cite{Ta}), whereas in
the hamiltonian case this relation makes no sense. Nevertheless,
Angenent and van der Vorst give an alternative description of the
index of a critical point $z=(u,v)$ in terms of the spectrum of an
integral operator associated with the matrix
\[
P(x)=\left(%
\begin{array}{cc}
  H_{uu}(x, z(x)) & H_{vu}(x, z(x)) \\
  H_{uv}(x, z(x)) & H_{vv}(x, z(x)) \\
\end{array}%
\right)
\]
(see \cite{AvdV2}, Section 3). We don't known if this variational
description of the index could be somehow used to obtain estimates
on the growth of the minmax sequences associated to the functional
$J$, as in \cite{BL}.
\end{remark}
\section{Variational formulation.}
In this section we establish the functional analytic framework
needed to study problem \eqref{pertsystem} from the variational
point of view, and we give the variational formulation for
\eqref{pertsystem}. \\
We begin with the spaces $\Theta ^r (\Omega )$, which are defined
in terms of the domains of fractional powers of the Laplacian in
$L^2(\Omega )$ with zero Dirichlet boundary conditions, i.e.
\[
-\Delta :H^2 (\Omega )\cap H_0^1(\Omega ) \subset L^2(\Omega)
\rightarrow L^2(\Omega)
\]
where $H^2 (\Omega ), H_0^1(\Omega )$ are the usual Sobolev
spaces; namely $\Theta ^r (\Omega )=D((-\Delta )^{r/2})$ for
$0\leq r\leq 2$, and the corresponding operator is denoted by
$A^r$
\[
A^r = (-\Delta )^{r/2} : \Theta ^r (\Omega ) \rightarrow L^2
(\Omega).
\]
The spaces $\Theta ^r$ are Hilbert spaces with inner product and
associated norm
\begin{eqnarray*}
(u,v)_{\Theta ^r} &=& \int _{\Omega } A^r u A^r v dx=((-\Delta )^{r/2} u,  %
(-\Delta )^{r/2} v)_{L^2},\\
\| u\|_{\Theta ^r} ^2 &=&\int _{\Omega } |A^r u|^2 dx =\| (-\Delta
)^{r/2} u \| _{L^2} ^2 ,
\end{eqnarray*}
see Lions and Magenes \cite{LM}. Let us fix in $H^1_0 (\Omega )$ a
system of orthogonal and $L^2$-normalized eigenfunctions $\varphi
_1, \varphi _2, \varphi _3,...$, of $-\Delta $, $\varphi _1 >0$,
corresponding to positive eigenvalues $\lambda _1 <\lambda _2 \leq
\lambda _3 \leq ...\uparrow +\infty $, counted with their
multiplicity. Then, writing
\[
u=\sum _{k=1} ^{\infty } \xi _k \varphi _k,\;\;\;\; \textrm{with }
\xi _k=\int_{\Omega}u \varphi _k dx,
\]
it is well known that
\begin{equation}
\label{defA_r}
 A^r u=(-\Delta )^{r/2} u=\sum _{k=1} ^{\infty } \lambda _k
 ^{r/2} \xi _k \varphi _k ,
\end{equation}
with domain
\begin{equation}
\label{defTeta_r}
 \Theta ^r (\Omega ) =D((-\Delta )^{r/2})=\{ \sum _{k=1} ^{\infty } \xi _k \varphi
 _k \in L^2 (\Omega ): \sum _{k=1} ^{\infty } \lambda _k
 ^r \xi _k ^2 <\infty \},
\end{equation}
if $r\geq 0$. Then we can identify $\Theta ^r (\Omega )$ with the
space
\begin{equation}\label{defTeta_r1}
\bar{\omega } ^r =\{ \xi = \{ \xi _k\} _{k=1} ^{\infty } : \sum
_{k=1} ^{\infty } \lambda _k ^r \xi _k ^2 <\infty \}, \;\;\;\;(\xi
, \eta )_r =\sum _{k=1} ^{\infty } \lambda _k ^r \xi _k \eta _k,
\end{equation}
and
\begin{equation}\label{defTeta_r2}
(u,v)_{\Theta ^r}= ((-\Delta )^{r/2} u,(-\Delta )^{r/2} v)_{L^2 }
=(\xi ,\eta )_r , \;\;\;\;\;\;\;\; \| u\| _{\Theta ^r} =|\xi | _r.
\end{equation}
The spaces $\Theta^r(\Omega)$, with $r<0$, can be
introduced as a representation of the dual spaces
$\Theta^r(\Omega)'$, using the Fourier characterization
\eqref{defTeta_r} of $\Theta^r(\Omega)$ (see \cite{HvdV}). The
motivation to introduce these spaces is to extend $A(\U )=\int
\nabla u\nabla v$ to functions $u$ and $v$ with different
regularity properties, and to define an appropriate functional
associated to \eqref{pertsystem}: this approach has been
introduced by Hulshof and van der Vorst in \cite{HvdV}, and by de
Figueiredo and Felmer in \cite{dFF}, hence we will be brief. Let
us first consider the quadratic part. Using the previous
notations, the quadratic form $A(\U )$ can also be written as
$A(\U )=\int \nabla u\nabla v =\sum _{k=1} ^{\infty } \lambda_k\xi
_k \eta _k$, where $u=\sum _{k=1} ^{\infty } \xi _k \varphi _k$
and $v=\sum _{k=1} ^{\infty } \eta _k \varphi _k$. Hence, if we
define the product Hilbert spaces
\begin{equation} \label{E_r}
E^r (\Omega ) =\Theta ^r (\Omega ) \times \Theta ^{2-r} (\Omega )
, \;\;\;\;\;\;\; 0<r<2,
\end{equation}
the quadratic form $A(\U )$ uniquely extends to a selfadjoint
bounded linear operator $L:E^r (\Omega )\rightarrow E^r(\Omega )$
as follows:
\begin{eqnarray} \nonumber
\nonumber \sum _{k=1} ^{\infty } \lambda _k \xi _k \eta _k &=&
\frac12 \sum _{k=1} ^{\infty } \lambda_k^r(\lambda_k^{1-r}\eta_k)
\xi _k + \frac12 \sum _{k=1}
^{\infty } \lambda_k^{2-r}(\lambda_k^{r-1}\xi_k) \eta _k\\
\nonumber &=&\frac12 ((-\Delta)^{1-r}v,u)_{\Theta^r}+\frac12
((-\Delta)^{r-1}v,u)_{\Theta^{2-r}} \\
\nonumber &=&\frac12 (L\U,\U)_{E^r}
\end{eqnarray}
where
\begin{equation} \label{defL}
L\U =((-\Delta )^{1-r} v, (-\Delta )^{r-1}u) \;\;\;\;\;
\U=(u,v)\in E^r(\Omega )
\end{equation}
(see \cite{HvdV}). Next we consider the eigenvalue problem
\[
L\U =\lambda \U \;\;\;\; \textrm{in } E^r(\Omega ).
\]
Using \eqref{defL} we can write equivalently
\begin{eqnarray*}
(-\Delta )^{1-r} v &=& \lambda u \\
(-\Delta )^{r-1} u &=& \lambda v
\end{eqnarray*}
which give directly
\[
v=\lambda ^2 v
\]
so that $\lambda = \pm 1$. The associated eigenvectors are
\begin{equation}\label{u+}
\U ^+ =(u,(-\Delta )^{r-1}u) \;\;\;\;\;\; \textrm{for }\lambda =1
\end{equation}
and
\begin{equation}\label{u-}
\U ^- =(u,-(-\Delta )^{r-1}u) \;\;\;\;\;\; \textrm{for }\lambda
=-1.
\end{equation}
We can define the eigenspaces
\begin{equation}\label{E+-}
E^{\pm } = \{ (u,\pm (-\Delta )^{r-1}u):u\in \Theta ^r (\Omega)
\};
\end{equation}
orthonormal bases consisting of eigenvectors of $E^{\pm}$ are
given by
\begin{equation}\label{Ebases}
\left\{ \mathbf{e}_k^{\pm}:= \frac{1}{\sqrt{2} } (\lambda _k^{-r/2} \varphi _k, \pm %
\lambda _k^{r/2-1} \varphi _k) \right\}_{k\in \mathbb{N} },
\end{equation}
and we have
\begin{equation}\label{orthogdecomp}
E^r(\Omega ) = E^+ \oplus E^- = \left\{ \U = \U ^+ + \U ^- , \U
^{\pm } \in E^{\pm } \right\}.
\end{equation}
We also find that, for $\U = \U ^+ + \U ^-$,
\[
A(\U )= \frac{1}{2} (L\U ,\U)_{E^r} = A(\U ^+) +A(\U^-),
\]
and
\[
A(\U ^+) - A(\U^-)=\frac{1}{2} \| \U \| ^2 _{E^r}.
\]
The derivative of $A(\U)$ defines a bilinear form on $E^r(\Omega)$
\begin{equation} \label{A'}
B(\U ,\mathbf{\Phi} ) =A'(\U )\mathbf{\Phi} = %
(L\U , \mathbf{\Phi} )_{E^r} ,\;\;\;\; \U ,\mathbf{\Phi} \in
E^r(\Omega).
\end{equation}
Next we define the Lagrangian $I(\U):E^r(\Omega)\rightarrow
\mathbb{R}$ associated to problem \eqref{pertsystem}. First of
all, we need the following Sobolev embedding theorem for
fractional order spaces (see \cite{LM} ):
\begin{theorem}\label{SobEmb}
If $0<2r<N$ the inclusions
\begin{equation}\label{sobolev}
\Theta ^r(\Omega) \hookrightarrow H^r(\Omega) \hookrightarrow
L^p(\Omega) \;\;\;\;\; \textrm{if }\;\; 1\leq p\leq
\frac{2N}{N-2r} <\infty
\end{equation}
are bounded; the second inclusion is compact if $1\leq p
<2N/(N-2r)$.\\
If $2r\geq N$, the inclusions are bounded and the second one is
compact for any $1\leq p< \infty$.
\end{theorem}
This theorem will allow us to define the Lagrangian associated to
problem \eqref{pertsystem} in a consistent way. An immediate
consequence of Theorem \ref{SobEmb} is that, by definition of
$E^r(\Omega)$,
\[
E^r(\Omega)\hookrightarrow L^{p+1}(\Omega)\times L^{q+1}(\Omega)
\]
if
\begin{equation}\label{E_rincl}
1\leq q+1 \leq
\dfrac{2N}{N-2r}, %
\;\;\;\;\; 1\leq p+1 \leq \dfrac{2N}{N+2r-4}
\end{equation}
for $N>2r$ and $N>4-2r$, that is,
\begin{equation}\label{range_r}
N\left[\frac12 -\frac1{q+1}\right] <r< 2-N\left[\frac12
-\frac1{p+1}\right].
\end{equation}
This embedding is compact if both inequalities bounding $p$ and
$q$ from above are strict. If $2r\geq N$, there is no restriction
on $p$, whereas if $4-2r\geq N$ there is no restriction on $q$.
Therefore, the Lagrangian $I$ associated to problem
\eqref{pertsystem} is well defined on $E^r(\Omega)$ if $p$ and $q$
satisfy inequality \eqref{E_rincl}, while we only require that
$0<r<2$: this restriction is motivated by the fact that we need
the compactness of the inclusion $E^r(\Omega) \hookrightarrow
L^2(\Omega)\times L^2(\Omega)$. The limiting values of $p$ and $q$
in \eqref{E_rincl} can be represented in the first quadrant of the
$(p,q)$-plane as a section of the well known critical hyperbola
\[
\frac 1{p+1} +\frac 1{q+1} =\frac {N-2}N
\]
which vanishes for $N\leq 2$.\\
Combining the extension $L$ of the quadratic form $A(\U)=\int
\nabla u\nabla v$ to $E^r(\Omega)$ defined in \eqref{defL} with
these inclusions, we can define the Lagrangian
\begin{eqnarray}\label{defI}
I(\U)&=&\frac 12 (L\U,\U)_{E^r} -\frac
1{q+1} \int_{\Omega} |u|^{q+1}dx-\frac 1{p+1} \int_{\Omega} |v|^{p+1}dx \\
\nonumber && -\int_{\Omega}kudx -\int_{\Omega}hvdx
\end{eqnarray}
associated to the perturbed system \eqref{pertsystem}, which is
well defined for $\U =(u,v) \in E^r(\Omega)$ if $p$, $q$ satisfy
\eqref{E_rincl} and $0<r<2$. We remark that critical points of
$I(\U)$ are classical solutions of problem \eqref{pertsystem}: see
for example \cite{HvdV}. Hence, to prove Theorem \ref{mainresult}
it suffices to show that $I(\U)$ has an unbounded sequence of
critical values. To do so, we require an estimate on the deviation
from symmetry of $I$ of the form
\begin{equation}\label{deviation}
|I(\U)-I(-\U)|\leq \beta (|I(\U)|^{1/\mu }+1)
\end{equation}
for $\U$ in $E^r(\Omega)$ and some $\beta >0$. Unfortunately $I$
does not satisfy \eqref{deviation}; however, it can be modified in
such a way that the new functional $J$ satisfy \eqref{deviation}
and large critical values of $J$ are also critical values of $I$.
\section{The symmetric case.}
In this section we consider the symmetric problem
\begin{equation}\label{symmsystem}
\left\{
\begin{array}{ll}
 -\Delta u= |v|^{p-1} v %
& \;\;\;\;\;x\in \Omega \\
& \\
 -\Delta v= |u|^{q-1} u %
& \;\;\;\;\;x\in \Omega \\
& \\
u=v=0 & \;\;\;\;\; x\in \partial \Omega
\end{array}
\right.
\end{equation}
that arises from \eqref{pertsystem} if $k(x) \equiv h(x) \equiv
0$. System \eqref{symmsystem} possesses a natural symmetry, which
guarantees the existence of infinitely many solutions. The aim of
this section is to exhibit these symmetrical critical values,
which will be used later on to construct the critical values of
the perturbed system \eqref{pertsystem}. The infinitely many
solutions of problem \eqref{symmsystem} can be found as critical
points of the corresponding functional $I$ by means of a version
of the symmetric Mountain Pass Theorem of Ambrosetti-Rabinowitz,
valid
for strongly indefinite functionals.\\
Let $E$ be  a Banach space with norm $\|\cdot\|$. Suppose that $E$
has a direct sum decomposition $E=E^1\oplus E^2$ with both $E^1$,
$E^2$ being infinite dimensional. Let $P^i$ denote the projections
from $E$ onto $E^i$. Assume $\{e^1_n\}$, $\{e^2_n\}$ are basis for
$E^1$ and $E^2$ respectively. Set
\begin{equation}\label{defX_n}
X_n=\langle e^1_1,...,e^1_n\rangle \oplus E^2, \hspace{20pt}%
X^k=E^1\oplus \langle e^2_1,...,e^2_k\rangle,
\end{equation}
and let $(X^k)^{\bot}$ denote the complement of $X^k$ in $E$. For
a functional $I\in\mathcal{C}^1(E,\mathbb{R})$ set $I_n:=I|_{X_n}$
the restriction of $I$ on $X_n$. Denote the upper and lower level
sets, respectively, by $I_a=\{z\in E:I(z)\geq a\}$, $I^b=\{z\in
E:I(z)\leq b\}$ and $I_a^b=I_a\cap I^b$. Then we have the
following theorem (see \cite{dFD}).
\begin{theorem}\label{simMP}
Let E as above and let $I\in \mathcal{C}^1(E, \mathbb{R})$ be even
with $I(0)=0$. In addition suppose, for each $k\in \mathbb{N}$,
the conditions below hold:
\begin{description}
\item [(I$_1$)] there is $R_k >0$ such that $I(\mathbf{z})\leq 0$
 for all $\mathbf{z}\in X^k$ with $\|\mathbf{z}\| \geq R_k$; \item
[(I$_2$)]there are $r_k>0$ and $a_k \rightarrow +\infty$ such that
$I(\mathbf{z})\geq a_k$ for all $\mathbf{z}\in (X^{k-1})^{\bot}$
with $\|\mathbf{z}\|=r_k$;
 \item [(I$_3$)] $I$ is bounded from
above on bounded sets of $X^m$;
 \item [(I$_4$)] $I$ satisfies the
$(PS)_c^{*}$ condition for any $c\geq 0$: that is, any sequence
$\{ \mathbf{z}_n \} \subset E$ such that $\mathbf{z}_n \in X_n$
for any $n \in \mathbb{N}$, $I(\mathbf{z}_n) \rightarrow c$ and
$I'_n(\mathbf{z}_n) \equiv \nabla(I| _{X_n})(\mathbf{z}_n)
\rightarrow 0$ as $n \rightarrow +\infty$ possesses a convergent
subsequence.
\end{description}
Then the functional $I$ possesses an unbounded sequence $\{c_k\}$
of critical values.
\end{theorem}
\begin{remark}
This theorem is a version of the Mountain Pass Theorem of
Ambrosetti and Rabinowitz for strongly indefinite symmetric
functionals, due to de Figueiredo and Ding (see \cite{dFD}). Other
versions of the same theorem are known, where the $(PS)_c^{*}$
condition is replaced by other variants, or by the usual (PS) (cf.
\cite{B}, \cite{BR}, \cite{D} and references therein).
\end{remark}
The sequence of critical values can be constructed by means of
certain Galerkin approximations (see \cite{BC}, \cite{BW},
\cite{dFD}), as we briefly recall. Using the previous notations,
set
\begin{equation}\label{defB_k}
B_k:= \left\{u\in X^k:\| \U\| \leq R_k \right\},
\end{equation}
the ball of radius $R_k$ in $X^k$,
\begin{equation}\label{defB_k^n}
B_k^n:= B_k \cap X_n = \left\{\U\in X^k \cap X_n:\| \U\| \leq R_k
\right\},
\end{equation}
and define the following sets of continuous maps
\begin{equation}\label{defGamma_k^n}
\Gamma_k^n:= \left\{h\in \mathcal{C} (B_k^n, X_n): h(-\U)=-h(\U),%
\hspace{3pt} h(\U)=\U \textrm{ on } \partial B_k^n \right\};
\end{equation}
finally define
\begin{equation}\label{c_k^n}
c_k^n:=\inf_{h\in \Gamma _k^n}{\sup_{\U \in B_k^n}{I(h(\U))}}.
\end{equation}
Then, it can be proved that for each $k\in \mathbb{N}$ fixed
(large enough, if necessary), the sequences $c_k^n$ converge to
critical values $c_k$ of the functional $I$, as $n$ tends to
$+\infty$; that is, the limits
\begin{equation}\label{defc_k}
c_k:= \lim_{n\rightarrow +\infty} c_k^n
\end{equation}
define critical values of the symmetric functional $I$. \\The
functional $I(\U)$ associated to the symmetric problem
\eqref{symmsystem} obviously satisfy the hypotheses of Theorem
\ref{simMP}, with $E=E^r(\Omega)$, $E^1=E^-$, $E^2=E^+$,
$e_j^1:=\mathbf{e}_j^-$, $e_j^2:=\mathbf{e}_j^+$, as we briefly
prove in the following. \bigskip \\
$\bullet \hspace{10 pt}$ As regards hypothesis (I$_1$), observe
that $\U \in X^k$ can be decomposed into the orthogonal sum
$\U=\U^+ +\U^-$ with $\U^- =(u_1, v_1) \in E^-$ and $\U^+ =(u_2
,v_2 )$ belonging to the finite dimensional space $\langle
\mathbf{e}_1^+, ...,\mathbf{e}_k^+ \rangle \equiv X^k \cap E^+$;
furthermore, by definition of the eigenvectors \eqref{Ebases},
$u_2$ and $v_2$ belong to the finite dimensional spaces
$E_k=\langle \varphi _1, ...\varphi _k \rangle $. By definition of
$I$, and recalling the embedding Theorem \ref{SobEmb}, we have,
for $\U \in X^k$, $\| \U \| _{E^r} =R$
\begin{eqnarray*}
I(\U) &=& \frac12 (L\U,
\U)_{E^r}-\frac1{q+1}\int_{\Omega}|u|^{q+1}dx
-\frac1{p+1}\int_{\Omega}|v|^{p+1}dx \\
&\leq& -\frac12 \| \U^-\|_{E^r}^2 + \frac12 \| \U^+\|_{E^r}^2
-c_q\int_{\Omega} \left(
|u_1|^{q+1} + |u_2|^{q+1} \right)dx\\
&&-c_p\int_{\Omega} \left( |v_1|^{p+1}
+ |v_2|^{p+1} \right) dx \\
&\leq& -\frac12 \| \U^-\|_{E^r}^2 -c_q\|u_1\|_{q+1}^{q+1}
-c_p\|v_1\|_{p+1}^{p+1} \\
&& + \frac12 \| \U^+\|_{E^r}^2%
-c_q \| u_2 \|_{\Theta ^r}^{q+1} \inf_{w \in E_k ,%
\hspace{5 pt} \|w\|_{\Theta ^r} =1}{\int_{\Omega} |w|^{q+1}dx} \\
&& - c_p \|v_2 \|_{\Theta ^{2-r}} ^{p+1} \inf_{w \in E_k,%
 \hspace{5 pt} \|w\|_{\Theta ^{2-r}} =1}{\int_{\Omega}
 |w|^{p+1}dx}\\
&\leq& \frac12 \| \U\|_{E^r}^2 - c_q(k,r)\| u_2
\|_{\Theta ^r}^{q+1} - c_p(k,r)\| v _2 \|_{\Theta ^{2-r}}^{p+1}\\
&\leq& \frac12 R^2 -c_{p,q}(k,r)R^{\min{(p+1,q+1)}}
\end{eqnarray*}
which tends to $-\infty$ as $R\rightarrow +\infty$. \bigskip \\
$\bullet \hspace{10 pt}$ The verification of hypothesis (I$_2$)
follows from the classical interpolation inequality in the $L^p$
spaces:
\begin{equation}\label{interpL^p}
\|f\|_{L^{s_0}}\leq \|f\|_{L^{s_1}}^{\alpha}\|f\|_{L^{s_2}}^{1-\alpha},%
\hspace{20 pt} \frac1{s_0} =\frac{\alpha}{s_1} + \frac{1-\alpha}{s_2}, %
0\leq \alpha \leq 1.
\end{equation}
Indeed, if $\U \in (X^{k-1})^{\bot}$, $(L\U, \U)_{E^r}=\| \U
\|_{E^r}^2$; combining the interpolation inequality
\eqref{interpL^p} with the Sobolev embedding \eqref{SobEmb} yields
\begin{eqnarray*}
I(\U) &=& \frac12 (L\U,
\U)_{E^r}-\frac1{q+1}\int_{\Omega}|u|^{q+1}dx
-\frac1{p+1}\int_{\Omega}|v|^{p+1}dx \\
&\geq& \frac 12 \| \U \|_{E^r}^2 -\left(\|u\|_2^{\alpha}%
\|u\|_{\frac{2N}{N-2r}}^{1-\alpha}\right) ^{q+1}%
-\left(\|v\|_2^{\beta} \|v\|_{\frac{2N}{N+2r-4}}^{1-\beta}\right)
^{p+1}
\end{eqnarray*}
where $\displaystyle \alpha = 1- \frac{N(q-1)}{2r(q+1)}$ and
$\displaystyle \beta = 1-\frac{N(p-1)}{(4-2r)(p+1)}$. Now, let us
observe that $(X^{k-1})^{\bot}=\langle \mathbf{e}_k^+,
\mathbf{e}_{k+1}^+, ...\rangle$: hence, if $\U=(u,v) \in
(X^{k-1})^{\bot}$ it is easy to verify, using definitions
\eqref{defTeta_r1}, \eqref{defTeta_r2} that
\begin{equation}\label{stimaortog_u}
\|u\|_2\leq \frac 1{\lambda _k^{r/2}}\|u\|_{\Theta^r},
\end{equation}
\begin{equation}\label{stimaortog_v}
\|v\|_2\leq \frac 1{\lambda _k^{1-r/2}}\|v\|_{\Theta^{2-r}}.
\end{equation}
Combining \eqref{stimaortog_u}, \eqref{stimaortog_v} with the
Sobolev embedding Theorem \ref{SobEmb} in the left hand side of
the previous inequality yields
\begin{eqnarray*}
I(\U)&\geq& \frac12
\|\U\|_{E^r}^2-\frac{C}{\lambda_k^{\alpha (q+1)r/2}}\|u\|_{\Theta^r}^{q+1}%
-\frac{C}{\lambda_k^{\beta
(p+1)(1-r/2)}}\|v\|_{\Theta^{2-r}}^{p+1}\\
&=&\frac12\|u\|_{\Theta^r}^2%
-C\left(\lambda_k^{-\frac{2r(q+1)-N(q-1)}{4(q+1)}}\|u\|_{\Theta^r}\right)^{q+1}\\
&&+\frac12\|v\|_{\Theta^{2-r}}^2%
-C\left(\lambda_k^{-\frac{(4-2r)(p+1)-N(p-1)}{4(p+1)}}\|v\|_{\Theta^{2-r}}\right)^{p+1}\\
&=&\|u\|_{\Theta^r}^2%
\left(\frac12-C\lambda_k^{-\frac{2r(q+1)-N(q-1)}{4}}\|u\|_{\Theta^r}^{q-1}\right)\\
&&+\|v\|_{\Theta^{2-r}}^2%
\left(\frac12-C\lambda_k^{-\frac{(4-2r)(p+1)-N(p-1)}{4}}\|v\|_{\Theta^{2-r}}^{p-1}\right).
\end{eqnarray*}
On the other hand, since $(X^{k-1})^{\bot}=\langle \mathbf{e}_k^+,
\mathbf{e}_{k+1}^+, ...\rangle\subseteq E^+$, by definition
\eqref{E+-} of the eigenspace $E^+$ if $\U=(u,v)\in
(X^{k-1})^{\bot}$, then $\U=(u,v)=(u, (-\Delta)^{r-1}u)$, and
\begin{equation}\label{relation_u_v}
\|v\|_{\Theta^{2-r}}=\|(-\Delta)^{r-1}u\|_{\Theta^{2-r}}=\|u\|_{\Theta^r}.
\end{equation}
Hence
\begin{eqnarray*}
I(\U)&\geq& \|u\|_{\Theta^r}^2%
\left(1-C\lambda_k^{-\frac{2r(q+1)-N(q-1)}{4}}\|u\|_{\Theta^r}^{q-1}\right. \\
&&\left.
-C\lambda_k^{-\frac{(4-2r)(p+1)-N(p-1)}{4}}\|u\|_{\Theta^r}^{p-1}\right)
\end{eqnarray*}
By \eqref{E_rincl}, the
exponents of $\lambda_k$ in the last expression are strictly
positive (recall that we choose $p,q$ below the critical
hyperbola); therefore, recalling that $\lambda_k \geq C\cdot
k^{2/N}$ for $k\rightarrow +\infty$, the
verification of (I$_2$) can be easily concluded. \bigskip \\
$\bullet \hspace{10 pt}$ Hypothesis (I$_3$) is clearly verified,
whereas the verification of $(PS)_c^{*}$ is standard, and follows
the one which will be given in the proof of ($\mathbf{3}$) of
Proposition \ref{Jproperties}, recalling that that $I(\U)\equiv
J(\U)$, so it is omitted here (see also \cite{HvdV} or \cite{D}).
\bigskip
\\Hence we can conclude that the symmetric problem
\eqref{symmsystem} possesses an unbounded sequence of critical
values, defined by \eqref{defc_k} and \eqref{c_k^n}.
\section{A modified functional.}
The aim of this section is to define a suitable modified
functional $J(\U)$, satisfying \eqref{deviation}, and whose
critical points are solutions of the original perturbed problem
\eqref{pertsystem}. We need first the following proposition.
\begin{proposition}\label{pr1}
There exists a constant $A$ depending on $\| h\| _{L^2(\Omega)}$,
$\| k\| _{L^2(\Omega)}$ such that if $I'(\U)\U =0$, then
\begin{equation}\label{est_crit}
\frac 1{q+1} \int_{\Omega} |u|^{q+1}dx+\frac 1{p+1} \int_{\Omega}
|v|^{p+1}dx \leq A\sqrt{I^2(\U) +1}.
\end{equation}
\end{proposition}
\begin{proof}[Proof of Proposition \ref{pr1}]. We follow the proof in
\cite{Ra}. Suppose that $I'(\U)\U =0$. Then, by simple estimates,
\begin{eqnarray*}
I(\U)&=& I(\U)-\frac 12 I'(\U)\U \\
&=&(\frac 12 -\frac 1{q+1})\int_{\Omega}|u|^{q+1}dx +
(\frac 12 -\frac 1{p+1})\int_{\Omega}|v|^{p+1}dx \\
&&-\frac 12 \int_{\Omega}h u dx -\frac 12 \int_{\Omega}k v dx \\
&\geq & C_1 \int_{\Omega}|u|^{q+1}dx +C_2 \int_{\Omega}|v|^{p+1}dx
-C_h \| u\|_2 -C_k \| v\|_2 \\
&\geq & C_3 \int_{\Omega}|u|^{q+1}dx +C_4 \int_{\Omega}|v|^{p+1}dx -C_5 \\
&\geq & C_6 \left\{ \int_{\Omega}|u|^{q+1}dx +
\int_{\Omega}|v|^{p+1}dx \right\} -C_7
\end{eqnarray*} where we
have used the following inequality: for any $\varepsilon>0 $ there
is a constant $C_\varepsilon >0$ such that
\[
\| f\|_2 \leq \varepsilon \|f \|_r^r +C_\varepsilon
\]
which is valid for any $f\in L^r(\Omega)$, $r>2$. \newline Hence
\eqref{est_crit} follows immediately.
\end{proof}
The idea underlying the construction of the modified functional
$J$ is, roughly speaking, to preserve the perturbation only where
$\int |u|^{q+1} + \int |v|^{p+1}$ is bounded from above by
$C|I(\U)|$, and to eliminate it where not.\newline To do so, let
$\chi \in \mathcal{C}^{\infty} (\mathbb{R}^+,\mathbb{R})$ be a
function satisfying $\chi (t)=1$ for $t\leq1$, $\chi (t)=0$ for
$t\geq 2$ and $-2<\chi '(t)<0$ for $1<t<2$. Set
\[
Q(\U)=Q(u,v)=2A\sqrt{I^2(\U)+1}
\]
and
\[
\psi(\U)=\psi(u,v)=\chi \left( \frac 1{Q(\U)}%
\left[ \frac 1{q+1} \int_{\Omega} |u|^{q+1}dx+\frac 1{p+1}
\int_{\Omega} |v|^{p+1}dx \right] \right).
\]
Note that if $\U$ is a critical point of $I$, then the argument of
$\chi$ lies in $[0,\frac 12 ]$ by Proposition \eqref{pr1} and
therefore $\psi(\U)=1$. Finally we set
\begin{eqnarray}\label{defJ}
J(\U)&=&\frac 12 (L\U,\U)_{E^r} -\frac
1{q+1} \int_{\Omega} |u|^{q+1}dx-\frac 1{p+1} \int_{\Omega} |v|^{p+1}dx \\
\nonumber && -\psi(\U) \left( \int_{\Omega}kudx +\int_{\Omega}hvdx
\right)
\end{eqnarray}
for $\U=(u,v)$ in $E^r(\Omega)$. It is easily seen that $J \in
\mathcal{C}^1 (E^r(\Omega),\mathbb{R})$; furthermore, if $\U$ is a
critical point of $I$, then $J(\U)=I(\U)$. the following
proposition contains the properties of $J$ which we need.
\begin{proposition}\label{Jproperties}
Let $f\in L^2(\Omega)$. Then
\begin{description}
  \item [(1)] There is a constant $\beta >0$ %
  depending on $\| h\|_2$ and $\| k\|_2$, such that
\begin{equation}\label{devJ}
|J(\U) -J(-\U)| \leq \beta \left(%
|J(\U)|^{\frac 1{q+1}}+ |J(\U)|^{\frac 1{p+1}} +1 \right)
\end{equation}
for $\U\in E^r(\Omega)$.
  \item [(2)] There is a constant $M_0>0$, depending on $\| h\|_2  $, $\|
  k\|_2$ such that \newline if $J(\U) \geq M_0$ and $J'(\U)=0$ %
  then $J(\U) =I(\U) $ and $I'(\U) =0$.
  \item [(3)] There is a constant $M_1\geq M_0$ such that for any $c>M_1$, \newline $J$
satisfies $\left( PS\right) _c$ and $(PS)_c^{*}$.
\end{description}
\end{proposition}
\begin{proof}[Proof of Proposition \ref{Jproperties}] We follow the proof
in \cite{Ra}, Proposition $10.16$.\bigskip \\
$\bullet \hspace{10 pt}$ To prove $(1)$, note first that if
$\U\notin \textrm{supp } \psi (\cdot) \cup \textrm{supp } \psi
(-\cdot)$, then $\psi(\U)=\psi(-\U)=0$ and $J(\U)=J(-\U)$, so that
\eqref{devJ} is valid. Hence, let us suppose that $\U\in
\textrm{supp } \psi (\cdot) \cup \textrm{supp } \psi (-\cdot)$. If
$\U\in \textrm{supp } \psi$, then
\begin{equation}\label{stimadev1}
\left| \int_{\Omega}kudx + \int_{\Omega} hvdx \right|\leq
\alpha_{q,p} \left( |I(\U)|^{\frac 1{q+1}} +|I(\U)|^{\frac
1{p+1}}+1\right),
\end{equation}
where $\alpha_{q,p}$ depends on $q$, $p$, $\| k\|_2$ and $\|
h\|_2$. Indeed, by Schwartz and H\"{o}lder inequalities and by
definition of $\psi(\U)$,
\begin{eqnarray*}
\left| \int_{\Omega}kudx + \int_{\Omega} hvdx \right| &\leq& %
\| k\|_2 \| u\|_2+\| h\|_2 \| v\|_2 %
\leq C\left( \| u\|_{q+1} +\| v\|_{p+1} \right) \\
&\leq& C \left( \int_{\Omega} |u|^{q+1} dx+%
\int_{\Omega}|v|^{p+1} dx \right)^{\frac 1{q+1}} +  \\
& & C \left( \int_{\Omega} |u|^{q+1} dx+
\int_{\Omega}|v|^{p+1} dx \right)^{\frac 1{p+1}}  \\
&\leq& C\left[ 4A\left(I^2(\U)+1\right)^{\frac 1{2(q+1)}} +%
4A\left(I^2(\U)+1\right)^{\frac 1{2(p+1)}} +1 \right]
\end{eqnarray*}
which implies directly \eqref{stimadev1}. Now, by definition,
\[
|J(\U)-J(-\U)|\leq (\psi(\U)+\psi(-\U))\left|\int_{\Omega}kudx
+\int_{\Omega} hvdx\right| ;
\]
combining this inequality with \eqref{stimadev1} yields
\begin{eqnarray*}
|J(\U)-J(-\U)|&\leq& \alpha_{q,p}%
(\psi(\U)+\psi(-\U))\left(|I(\U)|^{\frac
1{q+1}}+ |I(\U)|^{\frac 1{p+1}} +1 \right) \\
&\leq& c(\psi(\U)+\psi(-\U))%
\left(|J(\U)|^{\frac 1{q+1}}+ \left|\int_{\Omega}
kudx\right|^{\frac 1{q+1}}\right.\\
&& \left. +|J(\U)|^{\frac 1{p+1}}+ \left|\int_{\Omega}
hvdx\right|^{\frac 1{p+1}}+1 \right)\\
&\leq& 2c \left(|J(\U)|^{\frac 1{q+1}}+%
\left|\int_{\Omega}kudx\right|^{\frac 1{q+1}}%
 +|J(\U)|^{\frac 1{p+1}}+ \left|\int_{\Omega}
hvdx\right|^{\frac 1{p+1}}+1 \right).
\end{eqnarray*}
Since the exponents are smaller than 1, the $k$ and $h$ terms on
the right-hand side can be absorbed into the left-hand side
yielding
\eqref{devJ}. A similar estimate is valid for $\U\in \textrm{supp }\psi(-\cdot)$.\bigskip \\
$\bullet \hspace{10 pt}$ To prove $(\mathbf{2})$, it suffices to
show that if $M_0$ is large and $\U$ is a critical point of $J$
with $J(\U)\geq M_0$, then
\begin{equation}\label{stimaQ}
Q(\U)^{-1}\left( \frac1{q+1}\int_{\Omega}|u|^{q+1}dx+ %
\frac1{p+1}\int_{\Omega}|v|^{p+1}dx \right)<1;
\end{equation}
indeed, by definition of $\psi$, \eqref{stimaQ} implies $\psi
(\V)\equiv 1$ for $\V$ near $\U$. Hence $\psi '(\U)=0$, so
$J(\U)=I(\U)$, $J'(\U)=I'(\U)$ and $(\mathbf{2})$ follows.
Therefore we will prove that \eqref{stimaQ} holds. Let $\U=(u,v)$
and $\W=(w,z)$ be in $E^r(\Omega)$; hence, by definition of $J$,
\begin{eqnarray}\label{J'u}
J'(\U)\W&=&(L\U,\W)_{E^r}%
-\int_{\Omega}|u|^{q-1}u w dx -\int_{\Omega}|v|^{p-1}v z dx \\
&&\nonumber -\psi(\U)\left( \int_{\Omega}kw dx+\int_{\Omega}hz dx \right)%
-\psi'(\U)\W\left( \int_{\Omega}kudx+\int_{\Omega}hvdx \right),
\end{eqnarray}
where
\begin{eqnarray*}
\psi'(\U)\W&=&%
\chi'(\theta(\U))\theta'(\U)\W\\
&=&\chi'(\theta(\U))Q(\U)^{-2}%
\left\{ Q(\U)\left[ \int_{\Omega}|u|^{q-1}u w dx
+\int_{\Omega}|v|^{p-1}v z dx \right] \right. \\
&& \left. -(2A)^2 \theta(\U) I(\U)%
I'(\U)\W \right\}
\end{eqnarray*}
and
\[
\theta(\U)=Q(\U)^{-1}\left(%
\frac1{q+1}\int_{\Omega}|u|^{q+1}dx+\frac1{p+1}\int_{\Omega}|v|^{p+1}dx
\right).
\]
Regrouping terms in \eqref{J'u} yields
\begin{eqnarray}\label{J'uw}
J'(\U)\W&=&(1+T_1(\U))(L\U,\W)_{E^r}\\
\nonumber &&-(1+T_2(\U))\left(\int_{\Omega}|u|^{q-1}u w dx
+\int_{\Omega}|v|^{p-1}v zdx\right)\\
\nonumber &&-(\psi(\U)+T_1(\U))%
\left( \int_{\Omega}k w dx+\int_{\Omega}hz dx \right),
\end{eqnarray}
where
\begin{equation}\label{T_1}
T_1(\U)=\chi'(\theta(\U))(2A)^2\theta(\U)%
Q(\U)^{-2}I(\U)\left( \int_{\Omega}kudx+\int_{\Omega}hvdx \right)
\end{equation}
and
\begin{equation}\label{T_2}
T_2(\U)=T_1(\U)+\chi'(\theta(\U))%
Q(\U)^{-1}\left( \int_{\Omega}kudx+\int_{\Omega}hvdx \right).
\end{equation}
Let us now consider the term $J'(\U)\U$: from \eqref{J'uw} we have
\begin{eqnarray}\label{J'uu}
J'(\U)\U&=&(1+T_1(\U))(L\U,\U)_{E^r}\\
\nonumber &&-(1+T_2(\U))\left(\int_{\Omega}|u|^{q+1}dx
+\int_{\Omega}|v|^{p+1}dx\right)\\
\nonumber &&-(\psi(\U)+T_1(\U))%
\left( \int_{\Omega}k u dx+\int_{\Omega}h v dx \right);
\end{eqnarray}
therefore, if $\psi(\U)=1$ and $T_1(\U)=T_2(\U)=0$, we obtain
$J'(\U)\U=I'(\U)\U$ and $J(\U)=I(\U)$, so that \eqref{stimaQ}
follows from \eqref{est_crit}. Otherwise, consider
\begin{equation}\label{I-J'}
I(\U)-\frac1{2(1+T_1(\U))}J'(\U)\U
\end{equation}
and suppose that $\U$ is a critical point for $J$; since
$0\leq\psi(\U)\leq1$, if $T_1(\U)$ and $T_2(\U)$ are both small
enough, the calculation made in the proof of Proposition
\ref{pr1}, when carried out for \eqref{I-J'}, leads to
\eqref{est_crit} with $A$ replaced by a larger constant which is
smaller than $2A$: but then \eqref{stimaQ} holds. Therefore, it
suffices to show that $T_1(\U)$, $T_2(\U)\rightarrow 0$ as
$M_0\rightarrow \infty$.
\newline If $u\notin \textrm{supp }\psi$ then
$T_1(\U)=T_2(\U)=0$; hence we assume that $u\in \textrm{supp
}\psi$. Observe first that, by definition \eqref{T_1} of $T_1$ and
\eqref{stimadev1}
\begin{equation}\label{stimaT_1}
|T_1(\U)|\leq 4\alpha_{p.q}(|I(\U)|^{\frac1{q+1}}+%
|I(\U)|^{\frac1{p+1}}+1)|I(\U)|^{-1},
\end{equation}
where we have used the properties $|\chi '|<2$ and $\theta(\U)<2$
if $\U \in \textrm{supp }\psi$. Therefore, to conclude we need an
estimate relating $I(\U)$ and $J(\U)$ for $u\in \textrm{supp
}\psi$. By definition,
\[
I(\U)\geq J(\U)-\left| \int_{\Omega} kudx +\int_{\Omega} hvdx
\right|;
\]
thus, by \eqref{stimadev1},
\begin{equation}\label{IJ}
I(\U) +\alpha_{q,p} \left( |I(\U)|^{\frac 1{q+1}} +|I(\U)|^{\frac
1{p+1}}\right) \geq J(\U) -\alpha_{p,q} \geq M_0/2
\end{equation}
for $M_0$ large enough. If $I(\U)\leq 0$, estimate \eqref{IJ}
implies that
\[
\frac{\alpha_{p,q}^{(q+1)'}}{(q+1)'} + \frac{\alpha_{p,q}^{(p+1)'}}{(p+1)'}%
+\frac1{q+1}|I(\U)|+\frac1{p+1}|I(\U)|\geq M_0/2+ |I(\U)|,
\]
where $(q+1)'$, $(p+1)'$ are, respectively, the conjugate
exponents of $q+1$, $p+1$. But this is impossible for $p,q\geq 1$
and $M_0$ large enough: therefore, we can assume $I(\U)>0$. In
this case, \eqref{IJ} implies that $I(\U)\rightarrow +\infty$ as
$M_0 \rightarrow +\infty$, which shows, together with
\eqref{stimaT_1}, that $T_1(\U)\rightarrow 0$ as $M_0 \rightarrow
+\infty$. Analogous estimates yield $T_2(\U)\rightarrow 0$ as $M_0
\rightarrow +\infty$, and $(\mathbf{2})$ holds.\bigskip \\
$\bullet \hspace{10 pt}$ Let us now verify $(\mathbf{3})$. We have
to show that there is a constant $M_1\geq M_0$ such that for any
sequence $\{\U_n\}$ in $E^r(\Omega)$ satisfying
\begin{equation}\label{PSseq}
M_1<J(\U_n)<K \textrm{ for $n$ large}, \hspace{20pt}%
J'(\U_n)\rightarrow 0 \textrm{ as } n\rightarrow \infty
\end{equation}
has a convergent subsequence. The key point here is to prove that
such a sequence is necessarily bounded in $E^r(\Omega)$. Indeed,
by \eqref{J'uw} and \eqref{A'},
\[
J'(\U_n)=(1+T_1(\U_n))A'(\U_n) - \mathcal{H}(\U_n)
\]
where $\mathcal{H}$ is compact and $|T_1(\U _n)|\leq 1/2$ for
$M_1$ large enough; therefore, since $J'(\U_n)$ converges in
$(E^{r}(\Omega))'=E^{-r}(\Omega)$, the compactness of
$\mathcal{H}$ implies that a subsequence of $A'(\U_n)$ also
converges. In view of \eqref{A'}, we also have that $L\U_n$ and
$\U_n$ converge in $E^r(\Omega)$, because $L$ is invertible. To
prove that $\U_n$ is bounded we proceed as follows. By
\eqref{PSseq}, for any $\varepsilon >0$ there is a
$n(\varepsilon)$ such that for $n\geq n(\varepsilon)$
\begin{eqnarray*}
K+ \varepsilon \| \U_n \|_{E^r} &\geq&
J(\U_n)-\frac1{2(1+T_1(\U_n))}
J'(\U_n)\U_n \\
&=&\left(\frac{1+T_2(\U_n)}{1+T_1(\U_n)}\frac12-\frac1{q+1}\right)%
\int_{\Omega}|u_n|^{q+1}dx\\
&&+\left(\frac{1+T_2(\U_n)}{1+T_1(\U_n)}\frac12-\frac1{p+1}\right)
\int_{\Omega}|v_n|^{p+1}dx\\
&&+\left(\frac{\psi(\U_n)+T_1(\U_n)}{1+T_2(\U_n)}\frac12-\psi(\U_n)\right)%
\int_{\Omega}(ku_n+hv_n)dx;
\end{eqnarray*}
recalling that $T_1(\U), T_2(\U)\rightarrow 0$ as $M_1\rightarrow
+\infty$, we can choose $M_1$ sufficiently large such that the
coefficients of the integral terms $\int |u_n|^{q+1}$, $\int
|v_n|^{p+1}$ are strictly positive, that is (remember that $0\leq
\psi(\U_n)\leq 1$),
\begin{eqnarray*}
K+ \varepsilon \| \U_n \|_{E^r} &\geq& %
C_q\int_{\Omega}|u_n|^{q+1}dx+C_p\int_{\Omega}|v_n|^{p+1}dx%
-C\int_{\Omega}\left|ku_n+hv_n\right|dx \\
&\geq&C_q\int_{\Omega}|u_n|^{q+1}dx+C_p\int_{\Omega}|v_n|^{p+1}dx%
-C\left( \|k\|_2\|u_n\|_2+\|h\|_2\|v_n\|_2\right)\\
&\geq&C_q\int_{\Omega}|u_n|^{q+1}dx+C_p\int_{\Omega}|v_n|^{p+1}dx%
-C'\left( \|u_n\|_{q+1}+\|v_n\|_{p+1}\right)\\
&\geq&C'_q\int_{\Omega}|u_n|^{q+1}dx+C'_p\int_{\Omega}|v_n|^{p+1}dx-C''
\end{eqnarray*}
where the constants appearing in the previous inequalities depend
only on $M_1$, $q$, $p$ and not on $n$. Therefore we can conclude
that, for some new constants $K$, $\varepsilon >0$,
\begin{equation}\label{dimPS1}
K+\varepsilon \|\U_n\|_{E^r}\geq%
\int_{\Omega}|u_n|^{q+1}dx+\int_{\Omega}|v_n|^{p+1}dx.
\end{equation}
Decompose $\U_n=\U_n^{+} + \U_n^{-}$, where $\U_n^{+} \in E^{+}$,
$\U_n^{-} \in E^{-}$; writing $\U_n^{\pm}=(u_n^{\pm},v_n^{\pm})$,
we also have, by \eqref{u+}, \eqref{u-} and \eqref{orthogdecomp}
(the value of the constant $C$ can possibly change)
\begin{eqnarray*}
\|\U_n^{\pm}\|^2_{E^r}-\varepsilon
\|\U_n^{\pm}\|_{E^r} &\leq& \left| \left(%
L\U_n,\U_n^{\pm}\right)_{E^r}-\frac1{1+T_1(\U_n)}%
J'(\U_n)\U_n^{\pm}\right|\\
&&\\
&\leq&\frac{1+T_2(\U_n)}{1+T_1(\U_n)}\int_{\Omega}\left(%
|u_n|^q|u_n^{\pm}|+|v_n|^p|v_n^{\pm}|\right)dx\\
&&\\
&&+\frac{\psi(\U_n)+T_1(\U_n)}{1+T_1(\U_n)}%
\int_{\Omega}\left(|k||u_n^{\pm}|+|h||v_n^{\pm}|\right)dx\\
&&\\
&\leq&C\|u_n\|_{q+1}^q\|u_n^{\pm}\|_{q+1}+C\|v_n\|_{p+1}^p\|v_n^{\pm}\|_{p+1}\\
&&+C\left(\|u_n^{\pm}\|_2+|v_n^{\pm}\|_2\right)\\
&&\\
&\leq&C\|u_n\|_{q+1}^q\|u_n^{\pm}\|_{\Theta^r}+C\|v_n\|_{p+1}^p\|v_n^{\pm}\|_{\Theta^{2-r}}%
+C\|\U_n^{\pm}\|_{E^r}\\
&&\\
&\leq&C\left(\|u_n\|_{q+1}^q+\|v_n\|_{p+1}^p+1\right)\|\U_n^{\pm}\|_{E^r}.
\end{eqnarray*}
Dividing the first and the last expressions by
$\|\U_n^{\pm}\|_{E^r}$ we obtain
\begin{equation}\label{dimPS2}
\|\U_n^{\pm}\|_{E^r}-\varepsilon \leq
C\left(\|u_n\|_{q+1}^q+\|v_n\|_{p+1}^p+1\right).
\end{equation}
Combining \eqref{dimPS2} for $\U_n=\U_n^++\U_n^-$, together with
\eqref{dimPS1}, it follows that, possibly for some new constants,
\[
\|\U_n\|_{E^r}\leq C \{1+\{K+\varepsilon
\|\U_n\|_{E^r}\}^{q/(q+1)}+\{K+\varepsilon
\|\U_n\|_{E^r}\}^{p/(p+1)}\}
\]
which keeps $\|\U_n\|_{E^r}$ away from infinity. This implies that
the Palais-Smale condition is satisfied for $M_1$ sufficiently
large. The verification of $(PS)_c^{*}$ follows in the same way,
and thereby the proof is concluded.
\end{proof}
Property $(\mathbf{2})$ of Proposition $3.2$ guarantees that large
critical values of $J$ are also critical values of $I$; hence, in
what follows we shall seek large critical values of $J$.
\section{Minimax methods.}
The aim of this section is to construct suitable minimax sequences
which are strictly related to the existence of critical values of
the modified functional $J$, applying the method developed by
Rabinowitz to deal with perturbation from symmetry (see \cite{Ra})
to this case. The idea is to construct suitable minimax sequences
$d_k^n$, "perturbing" the ones defining the symmetric critical
values (in a sense that will be specified): comparison arguments
between the values of the two sequences will yield our
thesis.\\Observe first that, following the same lines as for the
verification of (I$_1$) of Theorem \ref{simMP}, it is not hard to
prove that for any $k\in \mathbb{N}$ there is a $R_k$ such that
$J(\U)\leq 0$ if $\U\in (B_k)^{\complement}$, where $B_k$ is the
sphere of radius $R_k$ in $X^k$ defined in \eqref{defB_k}. Hence
let us define the minimax sequences $c_k^n$ as in \eqref{c_k^n},
with $J(\U)$ instead of $I(\U)$, that is,
\[
c_k^n = \inf_{h\in \Gamma_k^n}{\sup_{\U \in B_k^n}{J(h(\U))}}.
\]
It is easy to verify that there is a sequence $b_k$, independent
on $n$, such that for any $k$ large enough
\begin{equation}\label{b_k}
c_k^n \leq b_k \hspace{20 pt} \textrm{ for any } n\in\mathbb{N}:
\end{equation}
indeed, by definition (since $\textrm{id} \in \Gamma _k^n$) and
applying Theorem \ref{SobEmb},
\begin{eqnarray*}
c_k^n &\leq& \sup_{\U\in B_k^n}{J(\U)} \leq \sup_{\U\in B_k}{J(\U)}\\
&\leq& \sup_{\U\in B_k}{\left[\frac12 \|\U ^{+}\|_{E^r}^2-\frac12
\|\U ^{-}\|_{E^r}^2 +C\|\U\|_2\right]}\\
&\leq&\sup_{\U\in B_k}{\left[\frac12 \|\U ^{+}\|_{E^r}^2-\frac12
\|\U ^{-}\|_{E^r}^2 +C\|\U\|_{E^r}\right]}\\
&\leq&\sup_{\U\in B_k}{\left[\frac12 \|\U
\|_{E^r}^2+C\|\U\|_2\right]}\leq CR_k^2
\end{eqnarray*}
for each $n\in\mathbb{N}$, and for $k \rightarrow +\infty$.
Following the idea in \cite{D}, \cite{BW}, \cite{dFD}, it is also
possible to prove that the sequence $c_k^n$ is bounded from below
by a sequence $a_k$ which is independent on $n$: that is, for any
$k$ large enough
\begin{equation}\label{a_k}
a_k\leq c_k^n \hspace{20 pt}\textrm{ for any } n\in \mathbb{N}:
\end{equation}
this fact, together with \eqref{b_k}, will be used to prove the
existence of the limit sequence $c_k$, as in the symmetric case.
To prove \eqref{a_k}, we need first the following version of the
Intersection Lemma:
\begin{lemma}\label{intersectionlemma}
Let us assume $B_k^n$, $\Gamma_k^n$ and $(X^{k-1})^{\bot}$ as
before; then, for any $h \in \Gamma_k^n$, $0<R<R_k$ there holds
\begin{equation}\label{intersection}
h(B_k^n)\cap \partial B_R \cap (X^{k-1})^{\bot} \neq \emptyset
\hspace{20 pt} \textrm{for any } n\in \mathbb{N}.
\end{equation}
\end{lemma}
\begin{proof}[Proof of Lemma \ref{intersectionlemma}] We follow the proof
given by Rabinowitz (Proposition $9.23$ in \cite{Ra}), hence we
will be brief. Let $\hat{O}_k^n:= \{ x \in B_k^n \mid h(x) \in B_R
\}$; since $h$ is odd, $0 \in \hat{O}_k^n$. Let $O_k^n$ denote the
component of $\hat{O}_k^n$ containing $0$. Since $B_k^n$ is
bounded, $O_k^n$ is a symmetric bounded neighborhood of $0$ in
$X^k \cap X_n$; therefore $\gamma (\partial O_k^n)=k+n$, where
$\gamma$ denotes the Krasnoleskii genus. We claim that
\begin{equation}\label{proofIL}
h(\partial O_k^n)\subset \partial B_R.
\end{equation}
Assuming \eqref{proofIL} for the moment. Set $W=\{ \mathbf{x}\in
B_k^n \mid h(\mathbf{x}) \in \partial B_R\}$; then \eqref{proofIL}
implies $\partial O_k^n \subset W$; hence, by the monotonicity
property of Krasnoleskii genus, $\gamma (W)=k+n$, so that $\gamma
(h(W)) \geq k+n$. Therefore, recalling that
$\textrm{codim}(X^{k-1})^{\bot} =k-1$, $h(W) \cap (X^{k-1})^{\bot}
\neq \emptyset$. On the other hand, by definition, $h(W) \subset
h(B_k^n) \cap \partial B_R$; consequently \eqref{intersection}
holds. \\
It remains to prove \eqref{proofIL}. Suppose $\mathbf{x} \in
\partial O_k^n$ and $h(\mathbf{x}) \in \stackrel{\circ }{B}_R$. If
$\mathbf{x} \in \stackrel{\circ }{B_k^n}$, there is a neighborhood
$N$ of $\mathbf{x}$ such that $h(N) \in \stackrel{\circ }{B}_R$:
but then $\mathbf{x} \notin \partial O_k^n$. Thus $\mathbf{x} \in
\partial B_k^n$, with $\partial$ relative to $X^n \cap X_k$; but
on $\partial B_k^n$, $h=id$. Consequently, if $\mathbf{x} \in
\partial B_k^n$ and $h(\mathbf{x}) \in \stackrel{\circ }{B}_R$,
$R>\| h(\mathbf{x})\| = \| \mathbf{x} \| =R_k$, contrary to the
hypothesis. Thus \eqref{proofIL} must hold.
\end{proof}
Applying the Intersection Lemma \ref{intersectionlemma} we are now
able to prove \eqref{a_k}. Let us fix $k$ large enough; then, for
any $n\in \mathbb{N}$, for any $h\in \Gamma_k^n$ and $0<R<R_k$
there is a $\W_n \in h(B_k^n) \cap \partial B_R\cap
(X^{k-1})^{\bot}$, so that by definition of $c_k^n$,
\begin{eqnarray*}
c_k^n &=& \inf_{h\in \Gamma_k^n}{\sup_{\U \in B_k^n}{J(h(\U))}}\\
&\geq& \inf_{h\in \Gamma_k^n}{J(\W_n)}\\
&\geq&\inf_{h\in \Gamma_k^n}{\sup_{0<R<R_k}{\inf_{\U \in \partial B_R \cap (X^{k-1})^{\bot}}{J(\U)}}}\\
&=&\sup_{0<R<R_k}{\inf_{\U \in \partial B_R \cap
(X^{k-1})^{\bot}}{J(\U)}};
\end{eqnarray*}
observe now that the last term of the previous inequality does not
depend on $n$, so that \eqref{a_k} is proved. Combining
\eqref{a_k} with \eqref{b_k} yields, for $k$ large enough,
\[
a_k\leq c_k^n \leq b_k \hspace{15 pt} \textrm{ for any } n \in
\mathbb{N}
\]
so that it is possible passing to the limit as $n$ tends to
$+\infty$ (up to a subsequence, if necessary), defining
\[
c_k=\lim_{n\rightarrow +\infty}{c_k^n}
\]
as in the symmetric case. This new minimax sequence $c_k$,
constructed for $J(\U)$, is not in general a sequence of critical
values for $J$, unless $k(x)\equiv h(x) \equiv 0$. \\ Let us now
construct new sequences  $d_k^n$, $d_k$ appropriately "perturbing"
$c_k^n$. First of all, we define a new sequence of sets
\begin{equation}\label{defU_k^n}
U_k^n:=\left\{ \U=t\mathbf{e}^+_{k+1}+\W \mid t\in [0,R_{k+1}], \W
\in B_{k+1}^n \cap X^k, \|u\| \leq R_{k+1} \right\};
\end{equation}
then define the new classes of functions
\begin{eqnarray}\label{Lambda_k^n}
\Lambda_k^n:&=&\left\{ H\in \mathcal{C}(U_k^n,X_n):H\mid_{B_k^n}
\in \Gamma_k^n , \hspace{3 pt} H(\U)=\U \right. \\
&& \nonumber \left. \textrm{ on } Q_k^n=(\partial B_{k+1}^n \cap
X^{k+1})\cup \left( (B_{R_{k+1}}\setminus B_{R_k})\cap X^k\right)
\right\}.
\end{eqnarray}
As one can easily observe, the new set $U_k^n$ is nothing that an
half of the sphere $B_{k+1}^n$, and the new class of function
$\Lambda_k^n$ are defined such that any $H\in \Lambda_k^n$,
suitably symmetrized, belongs also to $\Gamma_{k+1}^n$: combining
these facts with the estimate on the deviation from symmetry of
$J$, \eqref{devJ}, will be the key ingredient to obtain an upper
bound on the minimax sequences $c_k^n$ (and then also on $c_k$).
Now set
\begin{equation}\label{def_d_k^n}
d_k^n:=\inf_{H\in \Lambda_k^n}{\sup_{\U \in U_k^n}{J(H(\U))}}.
\end{equation}
Comparing the definition of $d_k^n$ with the one of $c_k^n$,
\eqref{c_k^n}, shows that $d_k^n \geq c_k^n$. Furthermore, we can
easily prove that $d_k^n$ is bounded from above independently on
$n$, as in \eqref{b_k}; indeed,
\begin{eqnarray*}
d_k^n &\leq& \sup_{\U\in U_k^n}{J(\U)} \leq \sup_{\U\in U_k}{J(\U)}\\
&\leq& \sup_{\U\in B_{k+1}}{\left[\frac12 \|\U
^{+}\|_{E^r}^2-\frac12
\|\U ^{-}\|_{E^r}^2 +C\|\U\|_2\right]}\\
&\leq&\sup_{\U\in B_{k+1}}{\left[\frac12 \|\U
^{+}\|_{E^r}^2-\frac12
\|\U ^{-}\|_{E^r}^2 +C\|\U\|_{E^r}\right]}\\
&\leq&\sup_{\U\in B_{k+1}}{\left[\frac12 \|\U
\|_{E^r}^2+C\|\U\|_2\right]}\leq CR_{k+1}^2
\end{eqnarray*}
for each $n\in\mathbb{N}$, and for $k \rightarrow +\infty$.
Therefore,
\[
a_k\leq c_k^n\leq d_k^n\leq \tilde{b}_k \hspace{10 pt} \textrm{
for any } n \in \mathbb{N}
\]
and it is possible to define (up to a subsequence)
\[
d_k =\lim_{n\rightarrow +\infty}{d_k^n};
\]
clearly, $d_k \geq c_k$. Furthermore, we have the following
fundamental proposition.
\begin{proposition}\label{prop_d_delta}
Assume $d_k>c_k\geq M_1$. For $\delta \in (0, d_k-c_k)$, define
\[
\Lambda_k^n(\delta) := \left\{ H\in \Lambda_k^n \mid J(H(\U))%
\leq c_k^n +\delta \textrm{ for } \U \in B_k^n \right\}
\]
and
\begin{equation}\label{defd_k^ndelta}
d_k^n (\delta):= \inf_{H\in \Lambda_k^n(\delta)}{\sup_{\U\in
U_k^n} J(H(\U))}.
\end{equation}
Then (eventually up to a subsequence) the limit
\begin{equation}\label{defd_kdelta}
d_k(\delta):= \lim_{n\rightarrow +\infty}{d_k^n(\delta)}
\end{equation}
exists for any $k\in \mathbb{N}$ large enough, and it is a
critical value of $J$.
\end{proposition}
\begin{proof}[Proof of Proposition \ref{prop_d_delta}]. The proof of
Proposition \ref{prop_d_delta} is based on the following, standard
"deformation lemma" (see, e.g., \cite{AR}).
\begin{lemma}\label{deformationlemma}
Let $E$ be a real Banach space, let $I\in \mathcal{C}^1(E,
\mathbb{R})$ and assume that $I$ satisfies $(PS)_c$. For $s\in
\mathbb{R}$ set $A_s=\{ u\in E \mid I(u)\leq s \}$. If $c$ is not
a critical value of $I$, given an $\bar{\varepsilon} >0$ there
exists an $\varepsilon \in (0, \bar{\varepsilon})$ and $\eta \in
\mathcal{C}([0,1]\times E,E)$ such that:
\begin{description}
\item [$1^{\circ}$] $\eta(t,u) =u$ for all $t\in [0,1]$, if $I(u)\notin %
[c-\bar{\varepsilon}, c+\bar{\varepsilon}]$ \item [$2^{\circ}$]
$\eta(1,A_{c+\varepsilon})\subset A_{c-\varepsilon}$.
\end{description}
\end{lemma}
The proof of Proposition \ref{prop_d_delta} follows the same lines
as in \cite{Ra}, adapted to this sort of Galerkin approximation
inspired to \cite{D}, \cite{BW} and others. If $d_k >c_k$, for any
$\delta \in (0, d_k -c_k)$ there is a $n_k \in \mathbb{N}$ (which
depends on $k, \delta$) such that
\[
0<\delta <d_k^n-c_k^n \hspace{20 pt} \textrm{ for any } n\geq n_k.
\]
Consider now, for $n>n_k$, $d_k^n(\delta)$ as defined in
\eqref{defd_k^ndelta}, and assume that it is not a critical value
of $J_n=J|_{X_n}$. Set $\bar{\varepsilon}= \frac12 (d_k^n -c_k^n
-\delta)>0$; then there exist $\varepsilon$ and $\eta$ as in the
deformation lemma \ref{deformationlemma}. Choose $H\in \Lambda
_k^n(\delta)$ such that
\begin{equation}\label{proofd_delta_1}
\max_{\U \in U_k^n}{J(H(\U))} \leq d_k^n(\delta) +\varepsilon.
\end{equation}
Consider $\eta(1, H(\cdot))$: clearly this function belongs to
$\mathcal{C}(U_k^n, X_n)$; if $\U \in Q_k^n$, $H(\U)=\U$ since $H
\in \Lambda_k^n$: therefore, $J(H(\U))=J(\U)\leq 0$ via the
definition of $R_k$ and $R_{k+1}$ (which do not depend on $n$).
Moreover, by the choice of $\bar{\varepsilon}$ and the assumption
$c_k >M_1>0$, $J(H(\U))=J(\U)\leq 0 < c_k^n +\bar{\varepsilon} <
d_k^n -\bar{\varepsilon}\leq c_k^n(\delta) - \bar{\varepsilon}$.
Hence, by $1^{\circ}$ of the deformation lemma
\ref{deformationlemma}, we have
\[
\eta (1, H(\U))=H(\U)= \U \hspace{15 pt} \textrm{ for } \U \in
Q_k^n.
\]
Further, since $H\in \Lambda_k^n (\delta)$, if $\U \in B_k^n$,
\[
J(H(\U))\leq c_k^n +\delta <d_k^n -\bar{\varepsilon}\leq %
d_k^n(\delta) -\bar{\varepsilon}
\]
by the choice of $\delta$ and $\bar{\varepsilon}$. Therefore,
again by $1^{\circ}$ of the deformation lemma
\ref{deformationlemma},
\[
\eta (1, H(\U))=H(\U) \hspace{15 pt} \textrm{ for } \U \in B_k^n,
\]
so that we can conclude that $\eta(1, H(\cdot)) \in
\Lambda_k^n(\delta)$. Thus, by definition of $d_k^n(\delta)$, we
get
\begin{equation}\label{proofd_delta_2}
d_k^n(\delta)\leq \max_{\U \in U_k^n}{J(\eta(1, H(\U)))}.
\end{equation}
On the other hand, \eqref{proofd_delta_1} and $2^{\circ}$ of Lemma
\ref{intersectionlemma} yields
\[
\max_{\U \in U_k^n}{J(\eta(1,H(\U)))}\leq d_k^n(\delta)
-\varepsilon,
\]
contrary to \eqref{proofd_delta_2}. Hence $d_k^n(\delta)$ is a
critical value of $J_n$. Now, let us apply the $(PS)_c^{*}$
condition, which is satisfied by $J(\U)$: indeed, we have just
proved, for any $k$ large enough, the existence of a sequence
$\{\Z_k^n\}\subset E^r$ such that for each $n\geq n_k$, $\Z_k^n
\in X_n$, $J'_n(\Z_k^n)=0$ and $J(\Z_k^n)=d_k^n(\delta)
\rightarrow d_k(\delta)>M_1$ as $n\rightarrow +\infty$ (the
existence of the limit $d_k(\delta)$, up to a subsequence, can be
easily proved, since $c_k^n\leq d_k^n(\delta)\leq d_k^n$). Hence,
${\Z_k^n}$ is a $(PS)^{*}_c$ sequence (with $c=d_k(\delta)>M_1$),
and by property $(\mathbf{I}_4)$ of Proposition \ref{Jproperties}
we can conclude that, along a subsequence, $\Z_k^n \rightarrow
\Z_k$ as $n\rightarrow +\infty$, with $J(\Z_k)=d_k(\delta)$ and
$J'(\Z_k)=0$. Hence $d_k(\delta) = \lim_{n\rightarrow
+\infty}{d_k^n(\delta)}$ (up to a subsequence) is a critical value
of $J(\U)$, and the proof is completed.
\end{proof}
On the basis of Proposition \ref{prop_d_delta}, to prove the
existence of infinitely many critical values for $J(\U)$ it
suffices to show that, up to a subsequence,
\[
d_k > c_k \geq M_1 \hspace{15 pt} \textrm{ for } k\in \mathbb{N},
\textrm{ and } c_k \rightarrow +\infty \textrm{ as } k\rightarrow
+\infty.
\]
This will be done in the following sections, estimating the growth
of $c_k$.
\section{A lower bound for $c_k$}
The aim of this section is to obtain an estimate from below on the
growth of the minimax sequence $c_k$. First of all, we will obtain
lower bounds for the minimax values $c_k^n$, and then also for the
sequence $c_k$: we recall that this sequence will in general
not consists of critical values of $J$, unless $k(x)\equiv h(x) \equiv 0$.\\
To estimate from below the growth of $c_k^n$ we follow the same
argument used to prove \eqref{a_k}, based on the Intersection
Lemma \ref{intersectionlemma} combined with the classical
interpolation inequality, in the same spirit of \cite{Ra}.
\begin{proposition}\label{lowerboundck}
Let $\frac1{q+1}+\frac1{p+1}>\frac{N-2}N$ and
\begin{equation}\label{cond_pq_1}
N\left[ \frac12 -\frac 1{q+1}\right] <r<2-N\left[ \frac12-
\frac1{p+1}\right] ;
\end{equation}
then there are $\gamma >0$  and $\tilde{k} \in \mathbb{N}$ such
that for all $k\geq \tilde{k}$,
\begin{equation}\label{lbc_k}
c_k\geq \gamma k^{2\alpha_r}
\end{equation}
where
\begin{eqnarray}\label{p_1q_1}
\alpha_r&=&\min(q_1,p_1),  \\
\nonumber \\
q_1&=&\frac{q+1}{q-1}\frac{r}{N}-\frac12,  \\
p_1&=&\frac{p+1}{p-1}\frac{2-r}{N} -\frac12 .
\end{eqnarray}
\end{proposition}
\begin{proof}[Proof of Proposition \ref{lowerboundck}] We remark here that the pair
$(p,q)$ lies below the critical hyperbola; for any fixed $(p,q)$,
the value of $r$, which identifies the space $E^r$, is not fixed,
but can be chosen in the range defined by \eqref{cond_pq_1} (see
Theorem \ref{SobEmb}). The aim of this proposition is to obtain a
lower bound for $c_k$ that depends only $r$, with $r$ unknown in
\eqref{cond_pq_1} . The "optimal" choice of $r$, in dependance of
$p,q$, will be a fundamental argument of the next Section.

First of all, we will prove a lower bound for the minimax
sequences $c_k^n$. Let $k \in \mathbb{N}$ be fixed. Let $h\in
\Gamma_k^n$ and $R<R_k$. By the Intersection Lemma
\ref{intersectionlemma}, for any $n \in \mathbb{N}$ there exists a
$\W _n \in h(B_k^n) \cap
\partial B_R \cap (X^{k-1})^{\bot}$, so that
\begin{equation}\label{proof_lb1}
\max_{\U \in B_k^n}J(h(\U))\geq J(\W_n) \geq %
\inf _{\U \in \partial B_R {\cap (X^{k-1})^{\bot}}}J(\U).
\end{equation}
\bigskip
Therefore, to obtain a lower bound for $c_k^n$ we have to estimate
$J(\U)$, where $\U \in \partial B_R \cap (X^{k-1})^{\bot}$ and
$0<R<R_k$. As remarked in the symmetric case (Section 3), if
$\U\in
\partial B_R \cap (X^{k-1})^{\bot}\subset E^+$ then
\[
\U =(u,v)=(u, (-\Delta)^{r-1}u)
\]
and
\begin{equation}\label{u_v}
\|v\|_{\Theta^{2-r}}=\|(-\Delta)^{r-1}u\|_{\Theta^{2-r}}=\|u\|_{\Theta^r}.
\end{equation}
Suitable combining the classical interpolation inequality
\eqref{interpL^p} with the Sobolev embedding Theorem \ref{SobEmb},
as in the symmetric case, yields the following inequalities (that
coincides with the classical Gagliardo-Nirenberg inequalities when
$r$ is an integer)
\begin{eqnarray}\label{proof_lb2}
\|u\|_{q+1} &\leq& C \|u\|_2 ^{\theta} \|u\|_{\Theta ^r}^{1-\theta} \hspace{10 pt}%
\textrm{ with } \hspace{10 pt} \theta =1- \frac{N}{r}\left(
\frac12 - \frac1{q+1}\right),\\
\label{proof_lb3} \|v\|_{q+1} &\leq& C \|v\|_2 ^{\zeta} \|v\|_{\Theta ^{2-r}}^{1-\zeta} \hspace{6 pt}%
\textrm{ with }\hspace{10 pt} \zeta =1- \frac{N}{2-r}\left(
\frac12 - \frac 1{p+1} \right).
\end{eqnarray}
Since $\U \in (X^{k-1})^{\bot}$, then $(L\U, \U)_{E^r}=\|
\U\|_{E^r}^2=2\|u\|_{\Theta^r}^2$ by \eqref{u_v}; combining
\eqref{proof_lb2}, \eqref{proof_lb3} and the estimates
\eqref{stimaortog_u}, \eqref{stimaortog_v} we obtain (we use the
same letter $C$ for different constants)
\begin{eqnarray}\label{proof_lb4}
\nonumber J(\U)&\geq& \frac12 \|\U\|_{E^r}^2%
-\frac1{q+1}\|u\|_{q+1}^{q+1}-\frac1{p+1}\|v\|_{p+1}^{p+1}-\|k\|_2\|u\|_2-\|h\|_2\|v\|_2 \\
\nonumber &\geq& \|u\|_{\Theta^r}^2 %
-C\|u\|_{q+1}^{q+1}-C\|v\|_{p+1}^{p+1}dx -C\\
\nonumber &\geq& \|u\|_{\Theta^r}^2 - C \left(\|u\|_2^{\theta}
\|u\|_{\Theta ^r}^{1-\theta}\right)^{q+1} - C \left( \|v\|_2
^{\zeta}
\|v\|_{\Theta ^{2-r}}^{1-\zeta} \right)^{p+1} -C \\
\nonumber &\geq& \|u\|_{\Theta^r}^2 - C \lambda_k^{-\frac{r}2
\theta (q+1)} \|u\|_{\Theta ^r}^{q+1} - C \lambda_k^{-\frac{2-r}2
\zeta
(p+1)} \|u\|_{\Theta ^r}^{p+1} -C \\
&\geq& \|u\|_{\Theta^r}^2 - C k^{-\frac{r}{N} \theta (q+1)}
\|u\|_{\Theta ^r}^{q+1} - C k^{-\frac{2-r}{N} \zeta (p+1)}
\|u\|_{\Theta ^r}^{p+1} -C
\end{eqnarray}
since $\lambda _k \geq C k^{2/N}$ for $k\rightarrow +\infty$,
where $\theta$, $\zeta$ satisfy conditions \eqref{proof_lb2},
\eqref{proof_lb3}. Inserting these values of $\theta$, $\zeta$ in
the right hand side of \eqref{proof_lb4} we obtain
\begin{equation}\label{proof_lb5}
J(\U) \geq \|u\|_{\Theta_r}^2 -
C\frac{\|u\|_{\Theta_r}^{q+1}}{k^{q_r}} -C
\frac{\|u\|_{\Theta_r}^{p+1}}{k^{p_r}} -C.
\end{equation}
where \begin{eqnarray*}
q_r&=&\frac{r}{N}(q+1)\left[1-\frac{N}{r}\left(\frac12
-\frac1{q+1}\right)\right]\\
p_r&=& \frac{2-r}{N}(p+1)\left[ 1-\frac{N}{2-r}\left(\frac12
-\frac1{p+1}\right)\right]
\end{eqnarray*}
To maximize the righthand side in \eqref{proof_lb5}, let us choose
\[
\|u\|_{\Theta^r}\asymp k^{\alpha},
\]
where $\alpha$ is unknown; then
\begin{eqnarray*}
J(\U) &\geq& \|u\|_{\Theta_r}^2 -
C\frac{\|u\|_{\Theta_r}^{q+1}}{k^{q_r}} -C
\frac{\|u\|_{\Theta_r}^{p+1}}{k^{p_r}} -C \\
&\asymp& k^{2\alpha} - k^{\alpha (q+1)-q_r} -k^{\alpha (p+1)-p_r}
-C.
\end{eqnarray*}
It is easy to verify that the optimal choice of $\alpha$ is
\[
\alpha=\alpha_r=\min\left(\frac{q_r}{q-1},\frac{p_r}{p-1}\right)
=\min\left(q_1,p_1\right);
\]
that is, for any $r$ in \eqref{cond_pq_1} and for any $\U \in
\partial B_R \cap (X^{k-1})^{\bot}$ with
$R=2\sqrt{\gamma}k^{\alpha_r}$,
\begin{eqnarray*}
J(\U) &\geq& 4\gamma k^{2\alpha_r} -
C\gamma^{q+1}k^{\frac{2q_r}{q-1}}
-C\gamma{p+1}k^{\frac{2p_r}{p-1}} -C\\
&\geq& \gamma k^{2\alpha_r}
\end{eqnarray*}
for $\gamma$ small enough. We remark that the condition
$R=2\sqrt{\gamma}k^{\alpha_r}<R_k$ is satisfied since $J(\U)<0$ if
$\|\U\|_{E^r}\geq R_k$, by definition of $R_k$, and this
contradicts the last inequality.

We are now ready to complete the proof. By \eqref{c_k^n} and
\eqref{proof_lb1}, for any $0<R<R_k$
\[
c_k^n=\inf_{h\in \Gamma_k^n}\sup_{\U \in B_k^n}J(h(\U))\geq%
\inf_{\U \in \partial B_R {\cap (X^{k-1})^{\bot}}}J(\U);
\]
choosing $R=R(k)=2\sqrt{\gamma}k^{\alpha_r}$ as before yields
\[
c_k^n \geq \gamma k^{2\alpha_r}
\]
for any $n\in \mathbb{N}$. Since the constants appearing in the
last estimate do not depend on $n$, as we have just remarked, we
can pass to the limit for $n \rightarrow \infty$, obtaining the
thesis.
\end{proof}
\section{Proof of Theorem \ref{mainresult}}
In this final section we shall complete the proof of Theorem
\ref{mainresult}. The idea of the proof is a reduction to the
absurd: basing on Proposition \ref{prop_d_delta}, we will assume
that $c_k =d_k$ for $k$ large, obtaining an upper bound on the
growth of $c_k$ which is in contrast with the lower bound proved
in Proposition \ref{lowerboundck}. Hence we will conclude that
$d_k >c_k$ for $k$ large, which yields the existence of an
unbounded sequence of critical values $d_k(\delta)$ for $J$, then
also for $I$. Therefore, we need first an estimate from above on
the growth of $c_k$ (under the assumptions that $c_k = d_k$).
\begin{proposition}\label{upperboundc_k}
If $c_k=d_k$ for all $k\geq k_1$, there exist two constants
$\alpha_1$, $\alpha_2 >0$ and $k_2 \geq k_1$ such that
\begin{equation}\label{ubc_k}
c_k \leq \alpha_1 k^{\frac{q+1}{q}} +\alpha_2k^{\frac{p+1}{p}}
\end{equation}
for all $k\geq k_2$.
\end{proposition}
\begin{proof}[Proof of Proposition \ref{upperboundc_k}]. We follow the
proof in \cite{Ra}. Let $k>k_1$; then, there is a sequence
$\{\varepsilon_n \}$ (depending on $k$) such that
\[
d_k^n \leq c_k^n + \varepsilon_n \hspace{10 pt} \textrm{and } \hspace{10 pt}%
\varepsilon_n \rightarrow 0  %
\hspace{5 pt} \textrm{ as } n\rightarrow +\infty.
\]
Let $\varepsilon >0$ and choose $n$ large enough such that
$\varepsilon_n <\varepsilon$. Then, choose $H \in \Lambda_k^n$
such that
\begin{equation}\label{proofMainr1}
\max_{\U \in U_k^n}{J(H(\U))}\leq d_k^n +\varepsilon \leq c_k^n
+2\varepsilon.
\end{equation}
Since $B_{k+1}^n = U_k^n \cup (-U_k^n)$, $H$ can be continuously
extended to $B_{k+1}^n$ as an odd function, still denoted with
$H$. Therefore, $H\in \Gamma_{k+1}^n$ and
\begin{equation}\label{proofMainr2}
c_{k+1}^n = \inf_{h\in \Gamma_{k+1}^n}{\max_{\U \in
B_{k+1}^n}{J(h(\U))}}\leq \max_{\U \in B_{k+1}^n}{J(H(\U))} =
J(H(\W _k^n))
\end{equation}
for some $\W_k^n \in B_{k+1}^n$. If $\W_k^n \in U_k^n$, by
\eqref{proofMainr1} and \eqref{proofMainr2},
\begin{equation}\label{proofMainr3}
c_{k+1}^n \leq J(H(\W _k^n)) \leq \max_{\U \in U_k^n}{J(H(\U))}
\leq c_k^n +2\varepsilon.
\end{equation}
If $\W_k^n \in -U_k^n$, by the oddness of $H$ and the estimate on
the deviation from symmetry \eqref{devJ}, we obtain
\eqref{proofMainr1} and \eqref{proofMainr2},
\begin{eqnarray*}
J(H(-\W_k^n))&=&J(-H(\W_k^n))\\
&\geq& J(H(\W_k^n)) -\beta \left(%
|J(H(\W_k^n))|^{\frac1{q+1}} +|J(H(\W_k^n))|^{\frac1{p+1}}
 +1\right).
\end{eqnarray*}
Since $c_k \rightarrow +\infty$ as $k+\infty$, \eqref{proofMainr2}
and the previous inequality imply that $J(-H(\W_k^n))>0$ for $n$
and $k$ large enough. Then, combining \eqref{proofMainr2}, the
estimate on deviation from symmetry \eqref{devJ} and the oddness
of $H$ yields
\begin{eqnarray}\label{proofMainr4}
\nonumber c_{k+1}^n &\leq& J(H(\W_k^n)) =J(-H(-\W_k^n))\\
\nonumber &\leq& J(H(-\W_k^n)) +\beta \left(%
|J(H(-\W_k^n))|^{\frac1{q+1}} +|J(H(-\W_k^n))|^{\frac1{p+1}}
 +1\right)\\
 &\leq& c_k^n +2\varepsilon +\beta\left(%
|c_k^n +2\varepsilon|^{\frac1{q+1}} +|c_k^n+
2\varepsilon|^{\frac1{p+1}} +1\right),
\end{eqnarray}
where we have used the fact that  if $\W_k^n \in -U_k^n$, then
$-\W_k^n \in U_k^n$ and $J(H(-\W_k^n))\leq c_k^n +2\varepsilon$,
by \eqref{proofMainr3}. Since $\varepsilon$ is arbitrary
(recalling that $\varepsilon_n \rightarrow 0$ as $n\rightarrow
+\infty$), \eqref{proofMainr3} and \eqref{proofMainr4} imply
\begin{equation}\label{proofMainr5}
c_{k+1}\leq c_k +\beta \left( c_k^{\frac1{q+1}}+c_k^{\frac1{p+1}}
+1 \right)
\end{equation}
for all $k$ large enough. Applying standard arguments, inequality
\eqref{proofMainr5} implies directly our thesis: see e.g.
\cite{Ra}.
\end{proof}
\begin{proof}[Proof of Theorem \ref{mainresult}] By Proposition
\ref{prop_d_delta}, if $d_k >c_k$ for $k$ large (up to a
subsequence, if necessary) there is a sequence of unbounded
critical values $d_k(\delta)$ for $J$, and then also for $I$, by
Proposition \ref{Jproperties}. Therefore, it suffices to show that
$d_k >c_k$ for $k$ large. On the contrary, let us assume that
$c_k=d_k$ as $k\rightarrow +\infty$. Then the last Proposition
\ref{upperboundc_k} assures the estimate from above \eqref{ubc_k}
on the growth of $c_k$, which depends only on $p,q$; on the other
hand, we have the estimate from below \eqref{lbc_k} proved in
Proposition \ref{lowerboundck}, which depends on $r\in
\left(N\left(\frac12 - \frac1{q+1}),2-N(\frac12 - \frac
1{p+1}\right)\right)$, for any pair $(p,q)$ below the critical
hyperbola. Theorem \ref{mainresult} will be proved if we choose
$r\in \left(N\left(\frac12 - \frac1{q+1}),2-N(\frac12 - \frac
1{p+1}\right)\right)$ such that
\begin{eqnarray}\label{fin_cond}
\nonumber \min{(2q_1, 2p_1)}&=& \min{ \left(
2\frac{q+1}{q-1}\frac{r}{N}-1,\;\;
 2\frac{p+1}{p-1} \frac{2-r}{N} -1 \right) } \\
 &>& \max{\left\{\frac{q+1}{q},\frac{p+1}{p}\right\}}.
\end{eqnarray}
Let us now consider the case
\begin{equation}\label{hp_q>p}
q\geq p
\end{equation}
that is,
\[
\max{\left\{\frac{q+1}{q},\frac{p+1}{p}\right\}}= \frac{p+1}{p}.
\]
Our aim now is to discuss the value of $\max{(2q_1, 2p_1)}$.\bigskip \\

Let us fix $(p,q)$ below the critical hyperbola $ \frac1{p+1}
+\frac1{q+1}=\frac{N-2}{N}$. By definition of $p_1$ and $q_1$,
\[
\min{(2q_1, 2p_1)}= 2q_1
\]
if and only if
\[
2\frac{q+1}{q-1}\frac{r}{N}-1 \leq 2\frac{p+1}{p-1}\frac{2-r}{N}
-1,
\]
that is, as one can easily verify, if and only if
\begin{equation}\label{maxq1p1}
r\leq \frac{(p+1)(q-1)}{pq-1}:=r_{p,q}.
\end{equation}
First of all, we observe that this limiting value of $r$ is
consistent with the condition \eqref{cond_pq_1} on $r$; indeed,
\[
N\left[\frac12 -\frac1{q+1}\right] < \frac{(p+1)(q-1)}{pq-1}%
<2-N\left[\frac12 -\frac1{p+1}\right]
\]
for any pair $(p,q)$ such that
\[
\frac1{p+1} +\frac1{q+1} > \frac{N-2}{N}
\]
as one can verify (we have used the following decomposition
$pq-1=(p+1)(q+1)-(p+1)-(q+1)$). We have now two possible choices
for $r$:
\begin{description}
\item [$\mathbf{(i)}$] if $r\leq r_{p,q}$, then
$\min{(2q_1,2p_1)}=2q_1$; \item[$\mathbf{(ii)}$] otherwise, if we
choose $r> r_{p,q}$, then $\min{(2q_1,2p_1)}=2p_1$.
\end{description}
Hence we discuss separately the two cases.
\begin{itemize}
\item [(i)] Consider first the choice
\[
N\left[ \frac12 -\frac1{q+1}\right] <r \leq r_{p,q} =%
\frac{(p+1)(q-1)}{pq-1}.
\]
In this case, recalling that we are assuming $q\geq p$,
\eqref{fin_cond} yields
\[
2q_1 = \frac{2r}{N}\frac{q+1}{q-1} -1  > \frac{p+1}{p},
\]
that is,
\begin{equation}\label{cond_r_2}
r> \frac{N}2\frac{q-1}{q+1}\frac{2p+1}{p} := r^{L}_{p,q}.
\end{equation}
Observe first that
$r^{L}_{p,q}>N[\frac12-\frac1{q+1}]=\frac{N}2\frac{q-1}{q+1}$ for
any $p,q>1$, so that \eqref{cond_r_2} is a condition effectively
stronger than \eqref{cond_pq_1}. Hence, condition \eqref{cond_r_2}
can be satisfied for certain $r$ if and only if
\[
r^L_{p,q}<r_{p,q},
\]
that is,
\[
\frac{N}2\frac{q-1}{q+1}\frac{2p+1}{p}<\frac{(p+1)(q-1)}{pq-1}.
\]
Recalling that $pq-1=(p+1)(q+1)-(p+1)-(q+1)$, the last inequality
can be written as
\[
2-\frac2{N}<\frac2{p+1}+\frac2{q+1}+\frac1{p}\left(\frac1{p+1}+\frac1{q+1}-1\right),
\]
that gives the following condition on $(p,q)$:
\begin{equation}\label{fin_cond_pq_1}
\frac1{p+1}+\frac1{q+1}+\frac{p+1}{p(q+1)}>\frac{2N-2}{N}.
\end{equation}
Therefore, for any pair $(p,q)$ verifying \eqref{fin_cond_pq_1} we
can choose $r\leq r_{p,q}$ such that conditions \eqref{cond_pq_1},
\eqref{cond_r_2} are satisfied. Condition \eqref{fin_cond_pq_1}
defines a new region in the $(p,q)$ plane which is contained in
the subcritical region delimited by the critical hyperbola.
 \item [(ii)] Consider now the other
possible choice of $r$,
\[
r_{p,q}=\frac{(p+1)(q-1)}{pq-1} <r< 2-N\left[ \frac12-
\frac1{p+1}\right].
\]
In this case, recalling that we are assuming $q\geq p$,
\eqref{fin_cond} yields
\[
2p_1 = 2\frac{2-r}{N}\frac{p+1}{p-1}-1> \frac{p+1}{p},
\]
that is,
\begin{equation}\label{cond_r_3}
r<r_{p,q}^{U}:=2-\frac{N}2\frac{p-1}{p+1}\frac{2p+1}{p}.
\end{equation}
Observe that $r_{p,q}^{U}<2-N(\frac12 -\frac1{p+1})$ for any
$p>0$, so that \eqref{cond_r_3} is a condition effectively
stronger than \eqref{cond_pq_1}. Hence, condition \eqref{cond_r_3}
can be satisfied for certain $r$ if and only if
\[
r_{p,q}<r_{p,q}^{U},
\]
that is,
\[
\frac{(p+1)(q-1)}{pq-1}<2-\frac{N}2\frac{p-1}{p+1}\frac{2p+1}{p}.
\]
This inequality is equivalent to
\[
\frac{(q+1)(p-1)}{pq-1}>\frac{N}2\frac{p-1}{p+1}\frac{2p+1}{p};
\]
using the decomposition $pq-1=(p+1)(q+1)-(p+1)-(q+1)$ the last
inequality yields
\[
2-\frac2{N}<\frac2{p+1}+\frac2{q+1}+\frac1{p}\left(\frac1{p+1}+\frac1{q+1}-1\right),
\]
that is equal to condition \eqref{fin_cond_pq_1} found in case
(i).
\end{itemize}

From cases (i) and (ii), we can conclude that there are values of
$r$ (satisfying \eqref{cond_pq_1}) such that condition
\eqref{fin_cond} holds if and only if $(p,q)$ verify
\eqref{fin_cond_pq_1} (assuming  $q\geq p$). By symmetry
arguments, we immediately conclude that if $p\geq q$,
\eqref{fin_cond} holds for values of $(p,q)$ satisfying the
corresponding condition
\begin{equation}\label{fin_cond_pq_3}
\frac1{p+1}+\frac1{q+1}+\frac{q+1}{q(p+1)} >\frac{2N-2}{N}.
\end{equation}
Combining \eqref{fin_cond_pq_1} with \eqref{fin_cond_pq_3} yields
the thesis.
\end{proof}
\begin{remark}
We note that for $p=q$ the limiting curve \eqref{condTh} obtained
in Theorem \ref{mainresult} assumes the same (limiting) value
independently obtained by  Struwe and Rabinowitz in \cite{St1},
\cite{Ra2} in the case of a single perturbed equation.
Furthermore, we observe that the subcritical region in the $(p,q)$
plane obtained in Theorem \ref{mainresult} do not include
supercritical values (in the sense of Sobolev embedding) of $p$
and $q$, so that one could ask if the variational setting
introduced in Section 2 is meaningful. Nevertheless, the optimal
choice of the exponent $r$ (that is, the optimal choice of the
space $E^r$) is obtained not for $r=1$ (the case $E^1=H_0^1\times
H_0^1$): that is, the regularity allowed for the functions $u,v$
is strictly related to the pair $(p,q)$, and choosing a priori the
space $H_0^1\times H_0^1$ would be too restrictive also for
perturbed systems with Sobolev-subcritical nonlinear terms.
\end{remark}
\vspace{10pt} \textbf{Acknowledgments.} The author would like to
thank B. Ruf for some useful discussions about this problem.

\end{document}